\documentclass[review]{elsarticle}

\makeatletter
\def\ps@pprintTitle{%
	\let\@oddhead\@empty
	\let\@evenhead\@empty
	\def\@oddfoot{\centerline{\thepage}}%
	\let\@evenfoot\@oddfoot}
\makeatother
\usepackage{lineno,hyperref}
\usepackage{amsmath}
\usepackage{amsfonts}
\usepackage{float}
\usepackage{longtable}
\usepackage{multirow}
\newcommand{\myol}[2][3]{{}\mkern#1.5mu\overline{\mkern-#1mu#2}\mkern 0.5mu}
\newenvironment{VarDescription}[1]%
  {\begin{list}{}{%
    \settowidth{\labelwidth}{\textbf{#1}}%
    \setlength{\leftmargin}{\labelwidth}\addtolength{\leftmargin}{\labelsep}}}%
  {\end{list}}

\modulolinenumbers[1]

\journal{ }

\begin{document}

\begin{frontmatter}

\title{Representative Days for Expansion Decisions in Power Systems}

%%% Group authors per affiliation:


\author[mymainaddressLuis]{\'Alvaro Garc\'ia-Cerezo}
\ead{Alvaro.Garcia29@alu.uclm.es}

\author[mymainaddressLuis]{Luis Baringo\corref{mycorrespondingauthor}}\cortext[mycorrespondingauthor]{Corresponding author}
\ead{Luis.Baringo@uclm.es}

\address[mymainaddressLuis]{Escuela T\'{e}cnica Superior de Ingenieros Industriales de Ciudad Real, Department of Electrical Engineering, Universidad de Castilla–-La Mancha, Campus Universitario s/n, 13071 Ciudad Real, Spain}

\begin{abstract}
Short-term uncertainty should be properly modeled when the expansion planning problem in a power system is analyzed.
Since the use of all available historical data may lead to intractability, clustering algorithms should be applied in order to reduce computer workload without renouncing accuracy representation of historical data.
In this paper, we propose a modified version of the traditional K-means method that seeks to attain the representation of maximum and minimum values of input data, namely, the electric load and the renewable production in several locations of an electric energy system.
The crucial role of depicting extreme values of these parameters lies in the fact that they can have a great impact on the expansion and operation decisions taken.
The proposed method is based on the traditional K-means algorithm that represents the correlation between electric load and wind-power production.
Chronology of historical data, which influences the performance of some technologies, is characterized though representative days, each one composed of 24 operating conditions.
A realistic case study based on the generation and transmission expansion planning of the IEEE 24-bus Reliability Test System is analyzed applying representative days and comparing the results obtained using the traditional K-means technique and the proposed method.
\end{abstract}


\end{frontmatter}

%\linenumbers

\section*{Notation}

The main notation used in this paper is stated below for quick reference, while other symbols are defined as needed throughout the text.
A subscript $r$/$h$ in the symbols below denotes their values in the $r$th representative day/$h$th hour.

\subsection*{Indices}
\begin{VarDescription}{$\hspace{17mm} $}
\item[$d$] Demands.
\item[$g$] Conventional generating units.
\item[$h$] Hours.
\item[$\ell$] Transmission lines.
\item[$n$] Nodes.
\item[$r$] Representative days.
\item[$s$] Storage facilities.
\item[$w$] Wind-power units.

\end{VarDescription}

\subsection*{Sets}
\begin{VarDescription}{$\hspace{17mm} $}

\item[$RE\left(\ell\right)$] Receiving-end node of transmission line $\ell$.
\item[$SE\left(\ell\right)$] Sending-end node of transmission line $\ell$.
\item[$\varOmega_{n}^{{\rm{D}}}$] Demands located at node $n$.
\item[$\varOmega_{n}^{{\rm{G}}}$] Conventional generating units located at node $n$.
\item[$\varOmega_{n}^{{\rm{S}}}$] Storage units located at node $n$.
\item[$\varOmega_{n}^{{\rm{W}}}$] Wind-power units located at node $n$.
\item[$\varOmega^{{\rm{G+}}}$] Candidate conventional generating units.
\item[$\varOmega^{{\rm{L+}}}$] Candidate transmission lines.
\item[$\varOmega^{{\rm{S+}}}$] Candidate storage units.
\item[$\varOmega^{{\rm{W+}}}$] Candidate wind-power units.
\end{VarDescription}

\subsection*{Parameters}
\begin{VarDescription}{$\hspace{17mm} $}
\item[$B_{\ell}$] Susceptance of transmission line $\ell$ [$\Omega^{-1}$].
\item[$C_{g}^{{\rm{G}}}$] Operation cost of conventional generating unit $g$ [\$/MWh].
\item[$C_{d}^{{\rm{LS}}}$] Load-shedding cost of demand $d$ [\$/MWh].
\item[$E_{srh_{0}}^{{\rm{S}}}$] Energy initially stored in storage facility $s$ [MWh].
\item[$\myol{E}_{s}^{{\rm{S}}}$] Maximum level of energy of storage facility $s$ [MWh].
\item[$F$] Large enough positive constant.
\item[$I_{g}^{{\rm{G}}}$] Investment cost of candidate conventional generating unit $g$ [\$/MW].
\item[$\tilde{I}_{g}^{\rm{G}}$] Annualized investment cost of candidate conventional generating unit $g$ [\$/MW].
\item[$\myol{I}^{\rm{G}}$] Investment budget for building candidate conventional generating units [\$].
\item[$I_{\ell}^{{\rm{L}}}$] Investment cost of candidate transmission line $\ell$ [\$].
\item[$\tilde{I}_{\ell}^{\rm{L}}$] Annualized investment cost of candidate transmission line $\ell$ [\$].
\item[$\myol{I}^{\rm{L}}$] Investment budget for building candidate transmission lines [\$].
\item[$I_{s}^{{\rm{S}}}$] Investment cost of candidate storage facility $s$ [\$].
\item[$\tilde{I}_{s}^{\rm{S}}$] Annualized investment cost of candidate storage facility $s$ [\$].
\item[$\myol{I}^{\rm{S}}$] Investment budget for building candidate storage facilities [\$].
\item[$I_{w}^{{\rm{W}}}$] Investment cost of candidate wind-power unit $w$ [\$/MW].
\item[$\tilde{I}_{w}^{\rm{W}}$] Annualized investment cost of candidate wind-power unit $w$ [\$/MW].
\item[$\myol{I}^{\rm{W}}$] Investment budget for building candidate wind-power units [\$].
\item[$\myol{M}^{\rm{S}}_{s}$] Maximum number of units that can be built of candidate storage facility $s$.
\item[$\myol{P}^{\rm{D}}_{d}$] Peak power consumption of demand $d$ [MW].
\item[$\myol{P}^{\rm{G}}_{g}$] Capacity of conventional generating unit $g$ [MW].
\item[$\myol{P}^{\rm{L}}_{\ell}$] Capacity of transmission line $\ell$ [MW].
\item[$\myol{P}^{\rm{S}}_{s}$] Charging and discharging power capacity of storage facility $s$ [MW].
\item[$\myol{P}^{\rm{W}}_{w}$] Capacity of wind-power unit $w$ [MW].
\item[$\alpha_{wrh}$] Capacity factor of wind-power unit $w$ [pu].
\item[$\beta_{drh}$] Demand factor of demand $d$ [pu].
\item[$\Delta \tau$] Duration of time steps [h].
\item[$\eta_{s}^{\rm{S^{C}}}$] Charging efficiency of storage facility $s$.
\item[$\eta_{s}^{\rm{S^{D}}}$] Discharging efficiency of storage facility $s$.
\item[$\sigma_{r}$] Weight of representative day $r$ [days].
\end{VarDescription}

\subsection*{Optimization Variables}
\begin{VarDescription}{$\hspace{17mm} $}
	\item[$e_{s r h}^{\rm{S}}$] Energy stored in storage facility $s$ [MWh].
	\item[$m_{s}^{\rm{S}}$] Number of units to be built of candidate storage facility $s$.
	\item[$p_{g r h}^{{\rm{G}}}$] Power produced by conventional generating unit $g$ [MW].
	\item[$\overline{p}_{g}^{{\rm{G}}}$] Capacity to be built of conventional generating unit $g$ [MW].
	\item[$p_{\ell r h}^{{\rm{L}}}$] Power flow through transmission line $\ell$ [MW].
	\item[$p_{d r h}^{{\rm{LS}}}$] Load shed of demand $d$ [MW].
	\item[$p_{s r h}^{\rm{S^{C}}}$] Charging power of storage facility $s$ [MW].
	\item[$p_{s r h}^{\rm{S^{D}}}$] Discharging power of storage facility $s$ [MW].
	\item[$p_{w r h}^{{\rm{W}}}$] Power produced by wind-power unit $w$ [MW].
	\item[$\overline{p}_{w}^{{\rm{W}}}$] Capacity to be built of wind-power unit $w$ [MW].
	\item[$x_{\ell}^{\rm{L}}$] Binary variable that is equal to 1 if candidate transmission line $\ell$ is built, and 0 otherwise.
	\item[$\theta_{n r h}$] Voltage angle at node $n$ [rad].
\end{VarDescription}
%

%\newpage

\section{Introduction} \label{introduction}

The Generation and Transmission Expansion Planning (G\&TEP) problem is solve to determine the new facilities that should be built in a power system in order to ensure the supply of the electric load in the future, since the time frame of this problem can comprise several decades.
It is motivated by the growth in peak loads, the penetration of renewable generating units and the aging of transmission facilities.

In most electricity markets, a central entity is in charge of taking expansion decisions of the transmission network, i.e., which transmission lines should be built. The aim of this system operator is to minimize investment and operation costs preventing load-shedding.
In addition, expansion decisions of generating units are taken by private investors, whose purpose is to maximize their economic profits along with minimizing their financial risk.
Nevertheless, optimal solution is not guaranteed accounting the G\&TEP problem as two independent problems. This is the reason why the perspective of the system operator is generally considered in technical literature when dealing with a G\&TEP problem.
It means that the central entity attains the optimal solution minimizing operation and investment cost, of both transmission facilities and generating units.
Once this is done, the system operator must provide indications about optimal expansion of generating units, with whom the government should design policy plans to promote the investment  in certain technologies or locations.

Historical data are generally used to model the performance of power systems, since the more realistic the input data of the G\&TEP problem are, the more accurate the solution of the problem will be in comparison with the future situation.

Regarding short-term uncertainty, electric load and renewable production are the historical data whose variability is more important.
On the one hand, electric load is characterized by a daily evolution pattern.
Since its variability depends on human habits, its progression can be accurately predicted using historical data.
On the other hand, the generation of electric energy through renewable sources depends on meteorological conditions.
For instance, the performance of wind turbines depends on wind speed as well as electric energy produced by solar panels and hydroelectric power stations relies on sunlight and rainfall, respectively.

Note that short-term uncertainty associated with renewable generating units increases the complexity of G\&TEP problems due to weather forecast can be poorly predicted in advance as opposed to the daily evolution of electric load.
Hence, the inclusion of storage units in electric energy systems is required in order to improve the penetration of renewable generating units.
Thus, energy can be discharged from storage units when it is needed and stored when there is an excess energy.
In addition, electric load and renewable generation are dependent magnitudes; for instance, low electric load generally coincides in time with high wind-power production.
The optimization model used to solve G\&TEP problems should properly represent this correlation between electric load and renewable production.

Solving a G\&TEP problem commonly involves the use of hourly data especially when we consider technologies which depend on the chronology; for instance, storage units.
Nevertheless, the optimization problem can be intractable because of using large amount of historical data as input data.
Thus, it is required to reduce the amount of historical data used in order to achieve a near optimal solution in a reasonable time. For this purpose, several techniques have been implemented in technical literature, such as load-duration curves and the K-means method.

Load-duration curves technique depicts short-term uncertainty of electric load through different levels arranged into blocks, within each one of them an electric load cumulative distribution function is built.
Subsequently, these functions are divided into several sectors, which are respectively associated with a different probability, and we calculate the average value inside each one of them obtaining different levels of electric load.
This technique can be expanded considering electric load and renewable production; for instance, load- and wind-duration curves in the case of accounting wind-turbines in the electric energy system under study.
In this case, the performance of arranging renewable production data is equal to the method previously explained for electric load data.
Besides, both magnitudes share the same blocks, in where all combinations of different levels of electric load and renewable production can take place.
These combinations, that can be used as input data of optimization problems, receive the name of system operating conditions.
The accuracy of the solution obtained using them relies on the number of blocks and levels selected in this method, being greater when bigger these numbers are.
However, a commitment should be reached between accuracy and computation workload.
This criterion can be extrapolated to the rest of methods used in technical literature.
Duration curves have been used in many references in the technical literature, e.g., considering net load duration curves \cite{Caramanis82,Wogrin13} or load- and wind-duration curves \cite{Baringo12,Montoya15}

The K-means technique applies algorithms of arranging data into groups, whose centroids are used with the purpose of representing the input data as well as reducing computer workload. The weight of each centroid is associated with the number of input data inside its group.
This method has the advantage that, in contrast to load- and wind-duration curves technique, it can consider different correlations of electric load and renewable production in several locations of the electric energy system under study.
The K-means method is used, for example, in \cite{Baringo14,Dominguez15}.

Duration curves and traditional K-means methods are compared in \cite{Baringo13}.
Their main issue of these two methods is that it is not possible to include units with inter temporal constraints such as storage units in the expansion problems.
To deal with this issue, \cite{Dehghan16} proposes using a representative day of each season, while \cite{Nogales16} and \cite{Wogrin14} consider a modified K-means method.
The main drawback of these methods is that they may not represent accurately extreme values of input data.
In case of using electric load and renewable production as input data, maximum and minimum values can have a great effect on the solution of the optimization problem.

Within this context, the contributions of this paper are threefold:
\begin{enumerate}
\item To propose a modified version of the traditional K-means method to achieve that system operating conditions obtained as output data of this technique properly represent maximum and minimum values of input data.
\item To use this new method to obtain representative days of electric load and wind-power production, each one composed of 24 operating conditions, in order to characterize the chronology of the historical data and thus allow the inclusion of technologies that depend on the chronology, such as storage units, in the formulation of expansion models.
\item To provide and analyze the results of a realistic case study with the purpose of checking if the proposed method reaches an improvement in the outcomes in comparison with the traditional K-means method.
\end{enumerate}

The remaining of this paper is organized as follows.
Section \ref{SecMeth} explains the methodology of the traditional K-means method and the proposed modified version of this technique.
Section \ref{SecForm} provides the formulation of the G\&TEP problem.
Section \ref{SecCS} displays the results of a case study, where a comparison among the outcomes obtained applying the clustering methods mentioned above is analysed.
Finally, Section \ref{SecConc} concludes the paper with some relevant remarks.
%

%Survey: \cite{Merrick16}

%\section{Limitations of Current Approaches} \label{SecPF}

\section{Methodology} \label{SecMeth}

The K-means method is a clustering algorithm which aim is to arrange data into groups called clusters according to similarities. 
On the one hand, the inputs of this algorithm are historical data of two physical processes, namely, the electric load and the wind-power production in several locations of an electric energy system.
On the other hand, the outputs of this technique are the cluster centroids along with the number of observations located at each cluster.
Note that cluster centroids, defined by the values of the two physical processes involved, represent the system operating conditions, which can be used as input data in the resolution of optimization problems (e.g., a long-term planning problem).

The K-means technique is useful when dealing with a significant amount of data in optimization problems due to the reduction of computer workload.
In order to ensure this, the users of this method are able to choose the $K$ number of operating conditions which is obtained.
However, it must be taken into account that a low number of operating conditions can mean that the representation of the input data may not be very accurate. In contrast, a high number of clusters can lead to intractability.

\subsection{Input data}

It is important to normalize the input data before applying the algorithm in case of working with electric load and wind-power production data, because it is common that the order of magnitude of the first one is greater than in the case of the second one.
If the input data are not normalized and the orders of magnitude of the two parameters analyzed are not similar, the results of the clustering method can be influenced by the weight of one of the parameters at the time of computing the quadratic distances between each original observation and each cluster centroid.

At this point, it is necessary to note that operating conditions cannot represent the chronology of the historical data.
Due to the penetration of renewable generating units in the electric energy systems, the fact of not modeling the chronology of the input data can cause a distortion among the results obtained and the reality.
Therefore, in this paper we use representative days, each one composed of 24 operating conditions, in order to characterize the chronology of the historical data.
This means that technologies which depend on the chronology, such as storage units, can be included in the formulation of the expansion model.

\begin{figure}
	\centering
	\includegraphics[width=11 cm]{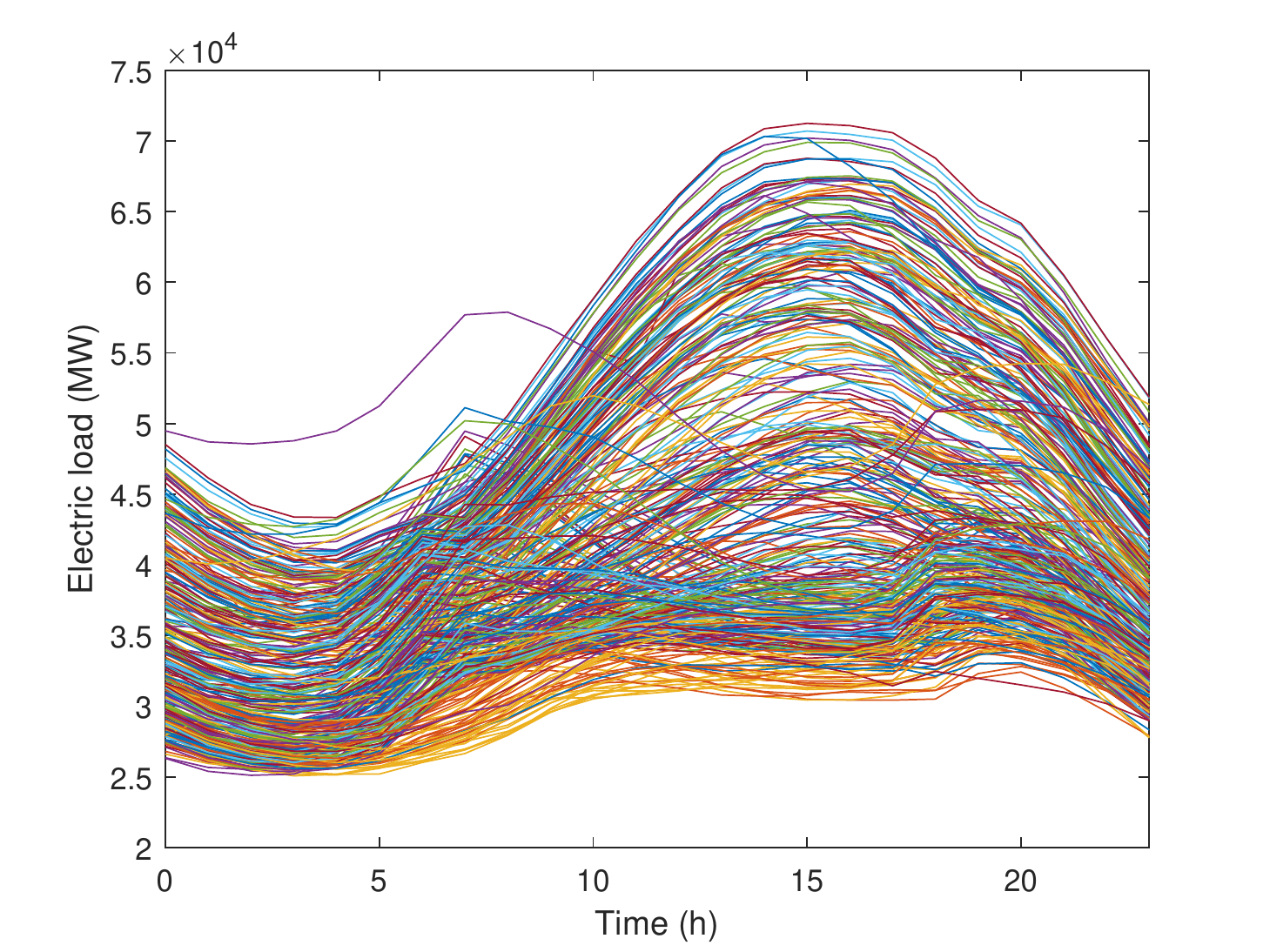}
	\caption{Daily evolution of electric load during a year.}
	\label{f301}
\end{figure}

\begin{figure}
	\centering
	\includegraphics[width=11 cm]{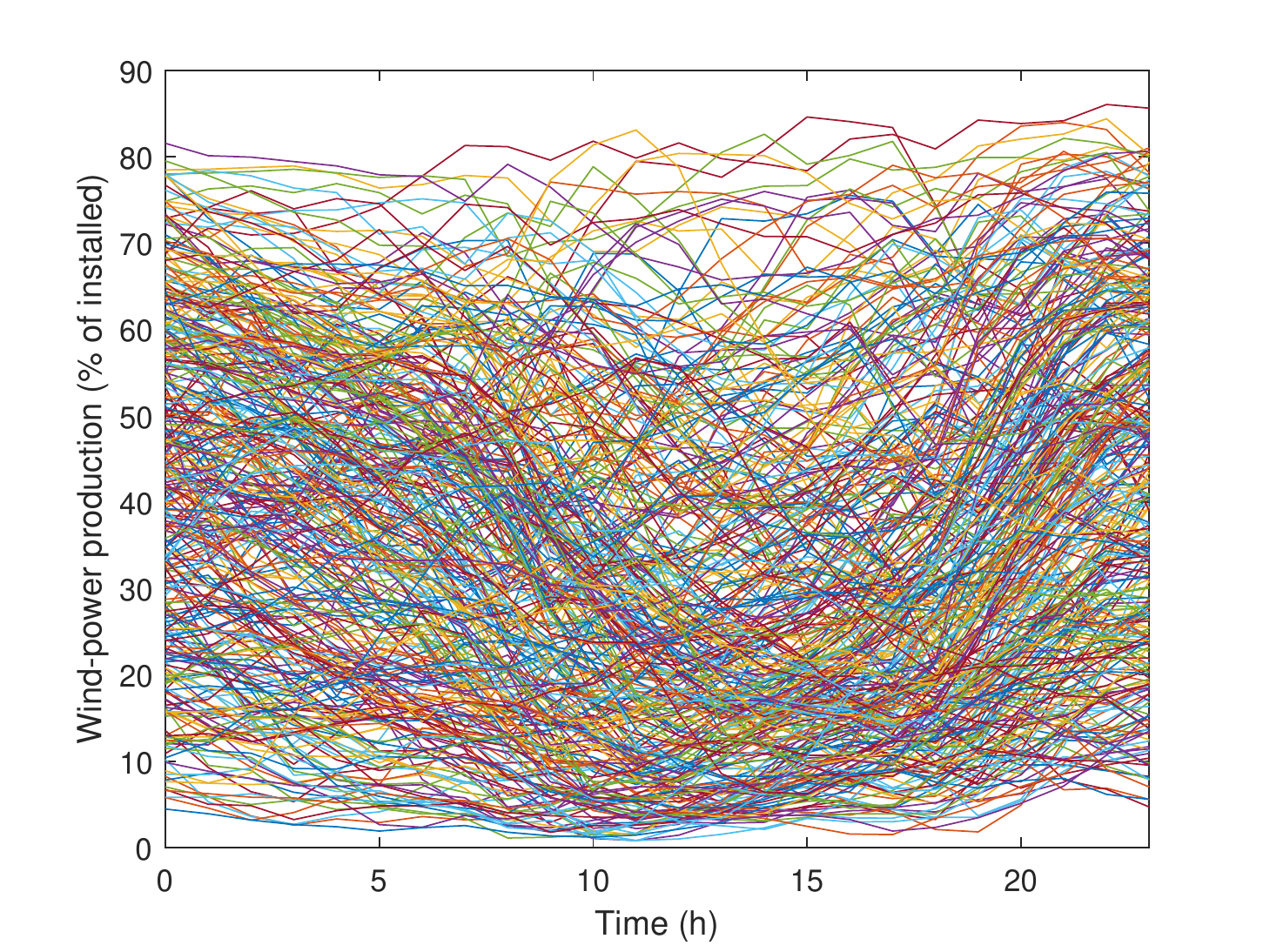}
	\caption{Daily evolution of wind-power production during a year.}
	\label{f302}
\end{figure}

We consider the historical data depicted in Figs. \ref{f301} and \ref{f302}, acquired from \cite{Energinet}, as input data of the algorithm.
Fig. \ref{f301} represents the daily evolution of electric load, while Fig. \ref{f302} illustrates the daily evolution of wind-power production, both during a year.
Note that, in this example, electric load units are MW, whilst wind-power production units are percent of installed. In Section \ref{SecCS}, the units considered for both parameters are MW.
A relevant aspect of Fig. \ref{f301} is that it displays a daily evolution pattern among different days of electric load data.
By contrast, Fig. \ref{f302} shows that the daily evolution of wind-power production does not follow any pattern.


\subsection{Traditional K-means algorithm}

The algorithm of the K-means method that has been used in technical literature, known from now on as traditional K-means method (TKM), is based on the following steps \cite{Baringo13}:
\begin{itemize}
	\item Step 1: Select the number of required clusters according to the needs of the problem.
	\item Step 2: Define the initial centroid of each cluster, e.g., randomly assigning a historical observation to each cluster.
	\item Step 3: Compute the quadratic distances between each original observation and each cluster centroid.
	\item Step 4: Allocate each historical observation to the closest cluster according to the distances calculated in Step 3.
	\item Step 5: Recalculate the cluster centroids using the historical observations allocated to each cluster.
\end{itemize}

\begin{figure}[H]
	\centering
	\includegraphics[width=8.5 cm]{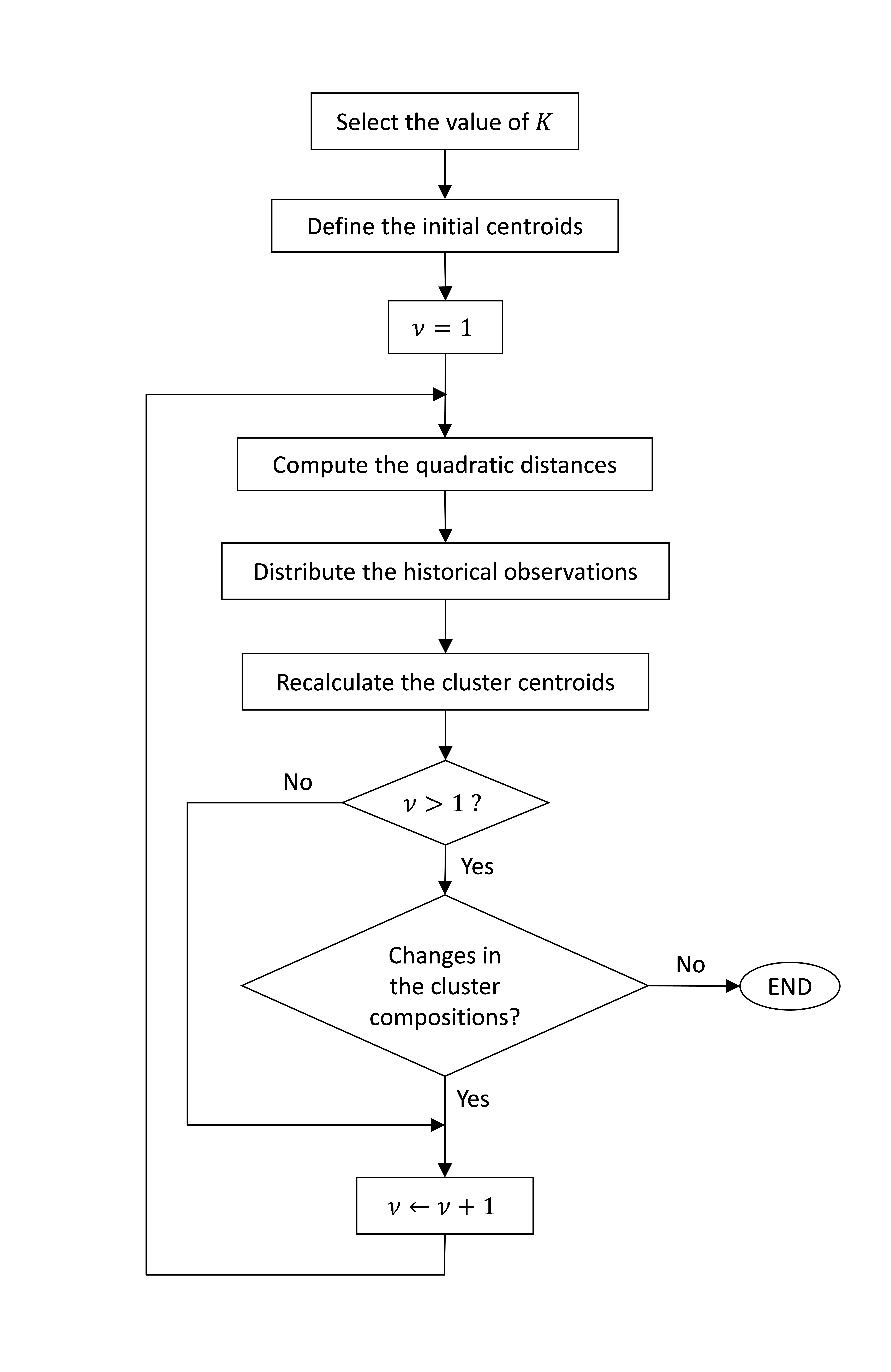}
	\caption{Flowchart of the traditional K-means method algorithm.}
	\label{flowchartTKM}
\end{figure}

Steps 3-5 are repeated iteratively until there are no changes in the cluster compositions between two consecutive iterations. Fig. \ref{flowchartTKM} illustrates the TKM algorithm.

In spite of the fact that the traditional K-means method presents advantages in comparison with other techniques (e.g., it is able to represent temporal and spatial correlations between uncertain parameters considered while duration curves technique cannot do it), it is not exempt of drawbacks.
The TKM sometimes does not adequately characterize the maximum and minimum values of the parameters analyzed.
This may constitute a problem, especially regarding the peak values, when we consider electric load and wind-power production as input data of the algorithm because their extreme values can have a great impact on the solution of the optimization problem.

In the case of the generation and transmission expansion planning (G\&TEP) problem, peak values of electric load can require the building of new generating units or new transmission lines to deliver the entire load of the electric energy system under study.
Not to mention that if the load cannot be completely supplied even then, the total costs will increase due to the load-shedding costs.
In addition, minimum values of electric load can also condition the solution of the optimization problem if the generating units have constraints linked to a minimum power produced greater than zero.
Moreover, extreme values of wind-power production can also influence the expansion and operation decisions taken.
Overall, maximum and minimum values of electric load and wind-power production can have an impact on the total costs, either by the investment costs associated with the expansion decisions made or by the operation costs related to the power produced by conventional generating units and load-shedding.

\subsection{Modified K-means algorithm}

To overcome these issues, we propose a new clustering method called modified K-means method (MKM), which tries to properly characterize the extreme values of the parameters considered, whose steps are presented below:
\begin{itemize}
	\item Step 1: Arrange the historical data into a $K_{1}$ number of clusters following the TKM.
	\item Step 2: Apply the same technique of clustering individually to the observations allocated to each cluster obtained in Step 1 arranging them into a $K_{2}$ number of clusters.
\end{itemize}

\begin{figure}[H]
	\centering
	\includegraphics[width=7 cm]{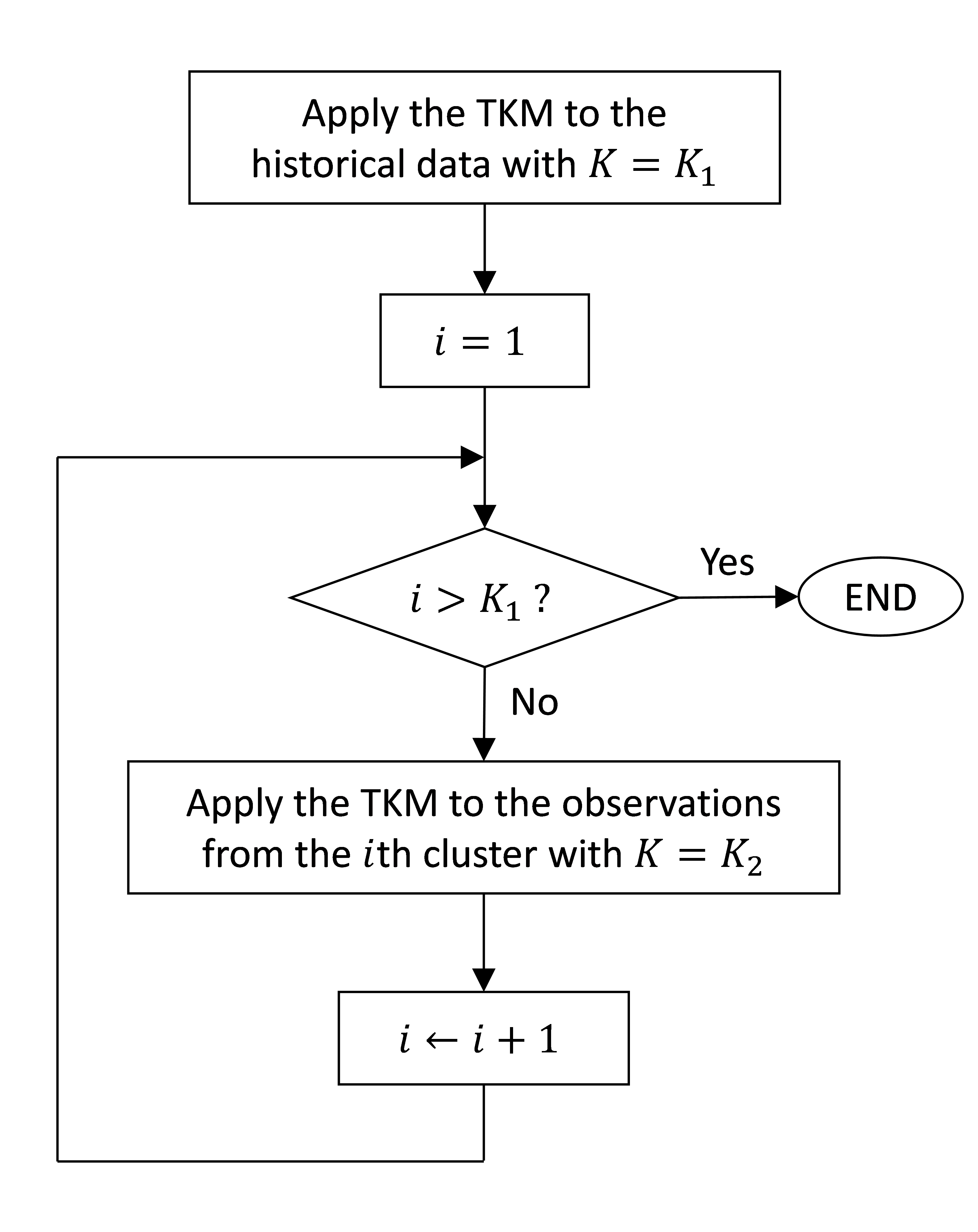}
	\caption{Flowchart of the modified K-means method algorithm.}
	\label{flowchartMKM}
\end{figure}

In other words, in Step 1 the MKM applies a first clustering to the historical data as it is customary in the technical literature, and then in Step 2 a second clustering is applied, but this time the input data are the observations of each cluster obtained in Step 1.
Therefore, Step 2 is repeated $K_{1}$ times until it has been applied to all the clusters acquired in the previous step.
The MKM algorithm is depicted in Fig. \ref{flowchartMKM}.

The number of operating conditions which are obtained as the output of this algorithm is equal to the product of $K_{1}$ and $K_{2}$.
For instance, a first clustering is applied organizing the input data into five clusters ($K_{1} = 5$). 
Then, the observations allocated to each cluster are considered as input data of a second clustering arranging them into two clusters ($K_{2} = 2$).
Thus, the number of operating conditions obtained at the end of the algorithm is 10.

Note that the MKM can only be applied if the number of observations located at each cluster after Step 1 is greater than or equal to the parameter $K_{2}$.
In addition, the parameter $K_{1}$ must be less than or equal to the number of input data considered in Step 1.
This last condition can be extrapolated to $K$ in the TKM.

Equation \eqref{eq301} defines the relation that must exist among the parameter $K$, associated with the traditional K-means method, and the parameters $K_{1}$ and $K_{2}$, linked to the modified K-means method, to make the results of both methods comparable.
\begin{equation}
K = K_{1}\cdot K_{2} \label{eq301}
\end{equation}

\subsection{Output data}

Since we use representative days in the case study described in Section \ref{SecCS}, we consider that the parameters $K$, $K_{1}$ and $K_{2}$ are associated with the number of representative days in their respective K-means methods, instead of the previous definitions that they have received in this paper.

\begin{figure}[H]
	\centering
	\includegraphics[width=11 cm]{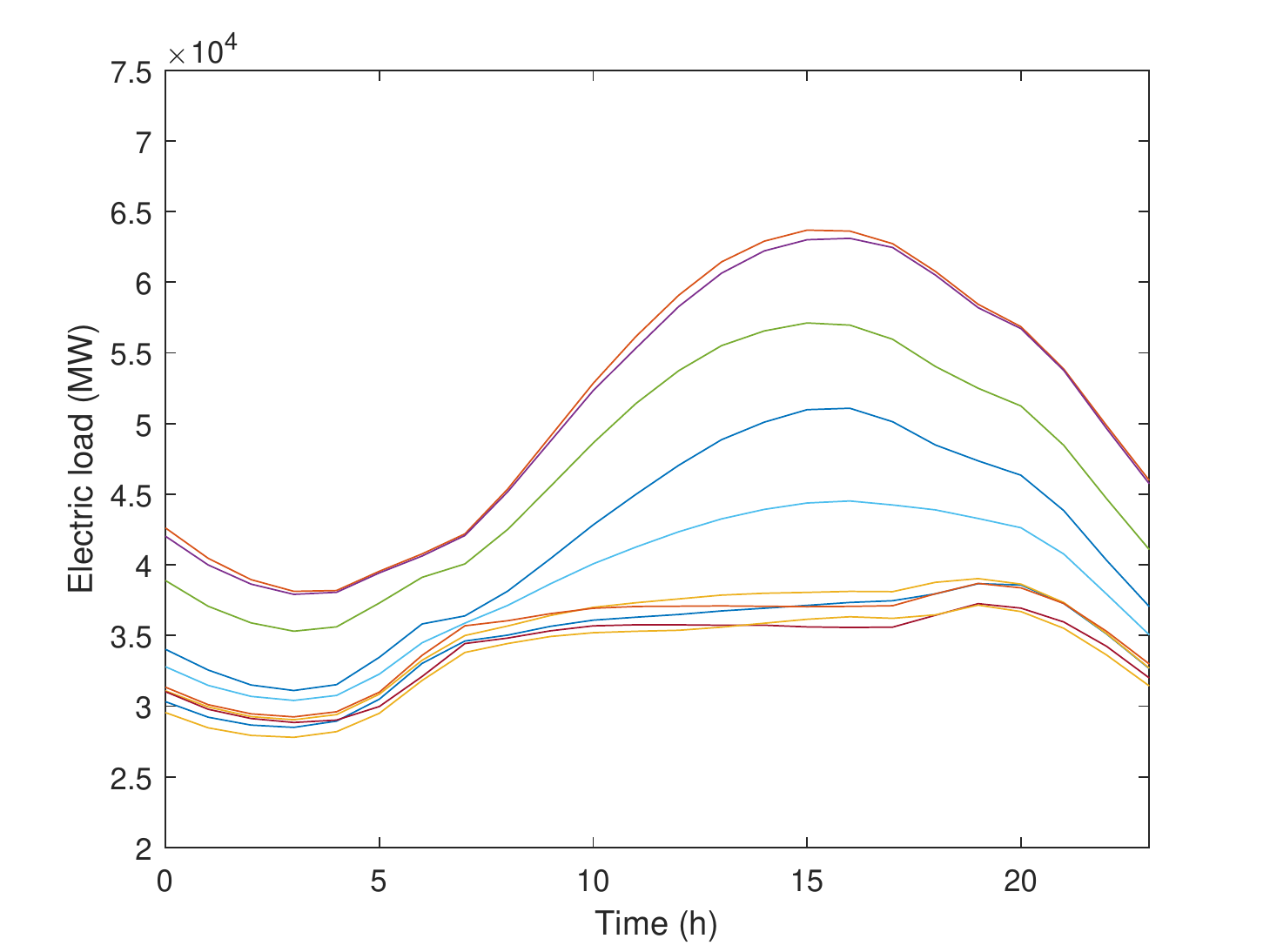}
	\caption{Representative days of electric load: traditional K-means method.}
	\label{f303}
\end{figure}

\begin{figure}[H]
	\centering
	\includegraphics[width=11 cm]{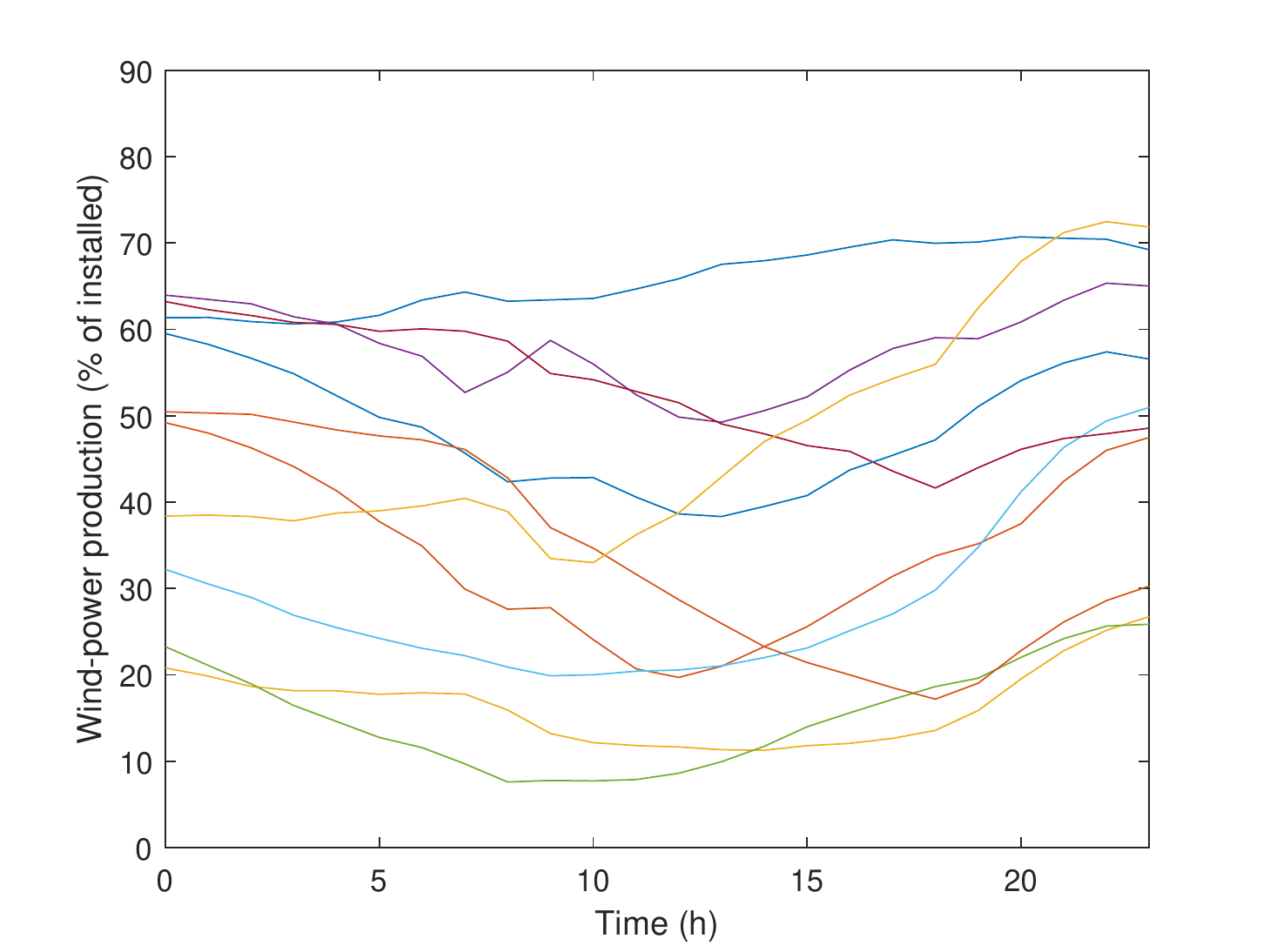}
	\caption{Representative days of wind-power production: traditional K-means method.}
	\label{f304}
\end{figure}

\begin{figure}[H]
	\centering
	\includegraphics[width=11 cm]{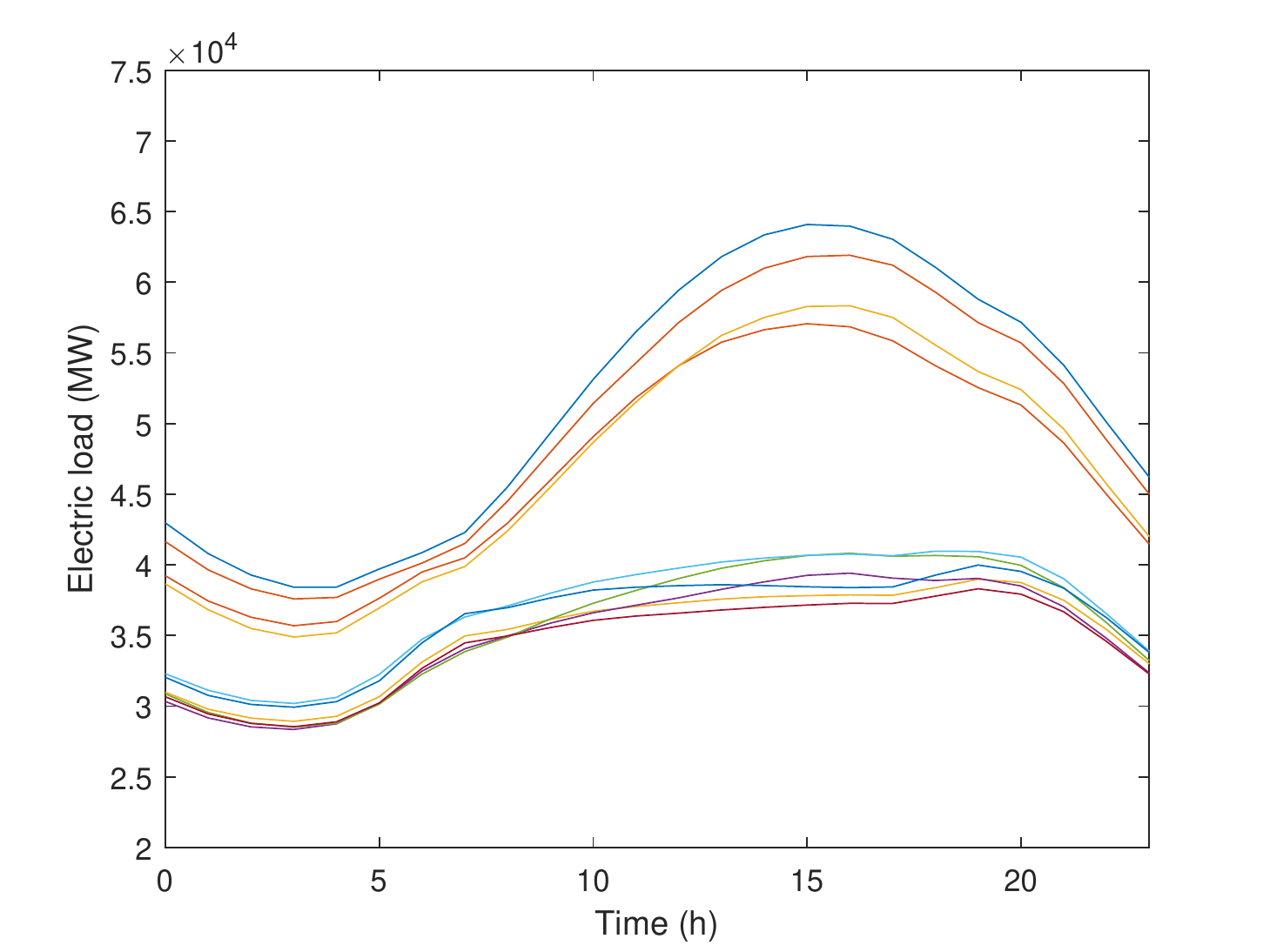}
	\caption{Representative days of electric load: modified K-means method.}
	\label{f305}
\end{figure}

\begin{figure}[H]
	\centering
	\includegraphics[width=11 cm]{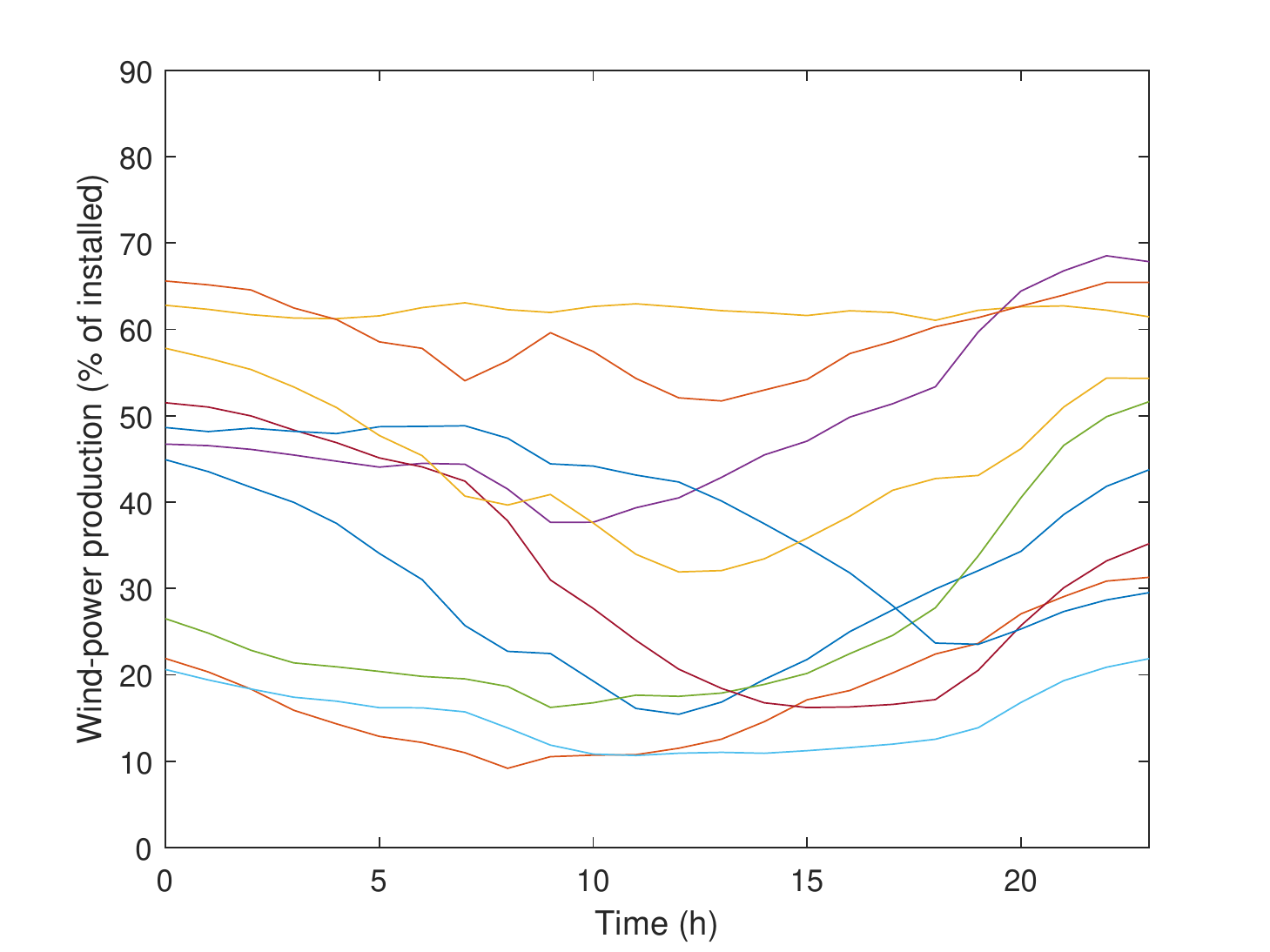}
	\caption{Representative days of wind-power production: modified K-means method.}
	\label{f306}
\end{figure}

The representative days of electric load and wind-power production obtained applying the traditional K-means method using $K = 10$ are illustrated in Figs. \ref{f303} and \ref{f304}, respectively.
Furthermore, Figs. \ref{f305} and \ref{f306} display the representative days obtained applying the modified K-means method using $K_{1} = 5$ and $K_{2} = 2$.
It is remarkable to mention that Fig. \ref{f305} shows more representative days of electric load in the areas of maximum and minimum values in comparison with Fig. \ref{f303}.


\section{Formulation} \label{SecForm}

The purpose of the G\&TEP problem is to minimize the operation costs along with the costs incurred in building new facilities (generating units, storage units, and transmission lines).
In this section, we provide the formulation of the G\&TEP problem considering a deterministic approach using the following mixed-integer nonlinear programming (MINLP) model:

\begin{subequations} \label{GandTEP}
\begin{align}
&{\rm{min}}_{\Phi}\quad\sum_{r}\sigma_{r}\sum_{h}\left(\sum_{g}C_{g}^{\rm{G}}p_{grh}^{\rm{G}} + \sum_{d}C_{d}^{\rm{LS}}p_{drh}^{\rm{LS}}\right) \nonumber \\
&+ \sum_{g\in\varOmega^{\rm{G+}}}\tilde{I}_{g}^{\rm{G}}\overline{p}_{g}^{\rm{G}} + \sum_{\ell\in\varOmega^{\rm{L+}}}\tilde{I}_{\ell}^{\rm{L}}x_{\ell}^{\rm{L}} + \sum_{s\in\varOmega^{\rm{S+}}}\tilde{I}_{s}^{\rm{S}}m_{s}^{\rm{S}}+ \sum_{w\in\varOmega^{\rm{W+}}}\tilde{I}_{w}^{\rm{W}}\overline{p}_{w}^{\rm{W}} \label{GandTEP1}
\end{align}
subject to
\allowdisplaybreaks
\begin{align}
&0\leq m_{s}^{\rm{S}}\leq \myol{M}_{s}^{\rm{S}}, \quad \forall s \in\varOmega^{\rm{S+}}, \label{GandTEP2}\\
&m_{s}^{\rm{S}} \in \mathbb{Z}, \quad \forall s \in \varOmega^{\rm{S+}}, \label{GandTEP3}\\
&0\leq \overline{p}_{g}^{\rm{G}} \leq \myol{P}_{g}^{\rm{G}}, \quad \forall g\in\varOmega^{\rm{G+}}, \label{GandTEP4}\\
&0\leq \overline{p}_{w}^{\rm{W}} \leq \myol{P}_{w}^{\rm{W}}, \quad \forall w\in\varOmega^{\rm{W+}}, \label{GandTEP5}\\
&x_{\ell}^{\rm{L}}\in \{ 0,1 \}, \quad \forall \ell \in \varOmega^{\rm{L+}}, \label{GandTEP6}\\
&\sum_{g\in\varOmega^{\rm{G+}}}I_{g}^{\rm{G}}\overline{p}_{g}^{\rm{G}}\leq \myol{I}^{\rm{G}}, \label{GandTEP7}\\
&\sum_{\ell\in\varOmega^{\rm{L+}}}I_{\ell}^{\rm{L}}x_{\ell}^{\rm{L}}\leq \myol{I}^{\rm{L}}, \label{GandTEP8}\\
&\sum_{s\in\varOmega^{\rm{S+}}}I_{s}^{\rm{S}}m_{s}^{\rm{S}}\leq \myol{I}^{\rm{S}}, \label{GandTEP9}\\
&\sum_{w\in\varOmega^{\rm{W+}}}I_{w}^{\rm{W}}\overline{p}_{w}^{\rm{W}}\leq \myol{I}^{\rm{W}}, \label{GandTEP10}\\
& \sum_{g\in\varOmega^{\rm{G}}_{n}}p_{grh}^{\rm{G}} +  \sum_{\ell|_{RE(\ell)=n}}p_{\ell rh}^{\rm{L}} + \sum_{s\in\varOmega^{\rm{S}}_{n}}p^{\rm{S^{D}}}_{srh} +
\sum_{w\in\varOmega^{\rm{W}}_{n}}p_{wrh}^{\rm{W}} = \nonumber\\
& \hspace{10mm}\sum_{d\in\varOmega^{\rm{D}}_{n}}(\beta_{drh}\myol{P}_{d}^{\rm{D}}-p_{drh}^{\rm{LS}}) + \sum_{\ell|_{SE(\ell)=n}}p_{\ell rh}^{\rm{L}} + \sum_{s\in\varOmega^{\rm{S}}_{n}}p^{\rm{S^{C}}}_{srh}, \quad \forall n, \forall r, \forall h, \label{GandTEP11}\\
&p^{\rm{L}}_{\ell rh} = B_{\ell}(\theta_{SE(\ell)rh}-\theta_{RE(\ell)rh}), \quad \forall \ell \setminus \ell \in \varOmega^{\rm{L+}} , \forall r, \forall h, \label{GandTEP12}\\
&p_{\ell rh}^{\rm{L}}=x_{\ell}^{\rm{L}}B_{\ell}(\theta_{SE(\ell)rh}-\theta_{RE(\ell)rh}), \quad \forall \ell \in \varOmega^{\rm{L+}}, \forall r, \forall h, \label{GandTEP13}\\
&-\myol{P}^{\rm{L}}_{\ell}\leq p_{\ell rh}^{\rm{L}} \leq \myol{P}^{\rm{L}}_{\ell}, \quad \forall \ell \setminus \ell \in \varOmega^{\rm{L+}} , \forall r, \forall h, \label{GandTEP14}\\
&-x_{\ell}^{\rm{L}}\myol{P}^{\rm{L}}_{\ell}\leq p_{\ell rh}^{\rm{L}}\leq x_{\ell}^{\rm{L}}\myol{P}^{\rm{L}}_{\ell}, \quad \forall \ell \in \varOmega^{\rm{L+}} , \forall r, \forall h, \label{GandTEP15}\\
&e_{srh}^{\rm{S}}=e_{srh-1}^{\rm{S}}+\left(  p_{srh}^{\rm{S^{C}}}\eta_{s}^{\rm{S^{C}}}-\frac{p_{srh}^{\rm{S^{D}}}}{\eta_{s}^{\rm{S^{D}}}}\right) \varDelta\tau, \quad \forall s, \forall r, \forall h > 1, \label{GandTEP16}\\
&e_{srh_{1}}^{\rm{S}}=E_{srh_{0}}^{\rm{S}}+\left(  p_{srh_{1}}^{\rm{S^{C}}}\eta_{s}^{\rm{S^{C}}}-\frac{p_{srh_{1}}^{\rm{S^{D}}}}{\eta_{s}^{\rm{S^{D}}}}\right) \varDelta\tau, \quad \forall s \setminus s \in \varOmega^{\rm{S+}}, \forall r, \label{GandTEP17}\\
&e_{srh_{1}}^{\rm{S}}=m_{s}^{\rm{S}}E_{srh_{0}}^{\rm{S}}+\left(  p_{srh_{1}}^{\rm{S^{C}}}\eta_{s}^{\rm{S^{C}}}-\frac{p_{srh_{1}}^{\rm{S^{D}}}}{\eta_{s}^{\rm{S^{D}}}}\right) \varDelta\tau, \quad \forall s \in \varOmega^{\rm{S+}}, \forall r, \label{GandTEP18}\\
&E_{srh_{0}}^{\rm{S}}\leq e_{srh_{24}}^{\rm{S}}, \quad \forall s \setminus s \in \varOmega^{\rm{S+}}, \forall r, \label{GandTEP19}\\
&m_{s}^{\rm{S}}E_{srh_{0}}^{\rm{S}}\leq e_{srh_{24}}^{\rm{S}}, \quad \forall s \in \varOmega^{\rm{S+}}, \forall r, \label{GandTEP20}\\
&0\leq e_{srh}^{\rm{S}}\leq \myol{E}_{s}^{\rm{S}}, \quad \forall s \setminus s \in \varOmega^{\rm{S+}}, \forall r, \forall h, \label{GandTEP21}\\
&0\leq e_{srh}^{\rm{S}}\leq m_{s}^{\rm{S}}\myol{E}_{s}^{\rm{S}}, \quad \forall s \in\varOmega^{\rm{S+}}, \forall r, \forall h, \label{GandTEP22}\\
&0\leq p_{grh}^{\rm{G}}\leq \myol{P}_{g}^{\rm{G}}, \quad \forall g \setminus g \in \varOmega^{\rm{G+}}, \forall r, \forall h, \label{GandTEP23}\\
&0\leq p_{grh}^{\rm{G}}\leq \overline{p}_{g}^{\rm{G}}, \quad \forall g\in\varOmega^{\rm{G+}}, \forall r, \forall h, \label{GandTEP24}\\
&0\leq p_{drh}^{\rm{LS}}\leq \beta_{drh}\myol{P}_{d}^{\rm{D}}, \quad \forall d, \forall r, \forall h, \label{GandTEP25}\\
&0\leq p_{srh}^{\rm{S^{C}}}\leq \myol{P}_{s}^{\rm{S}}, \quad \forall s \setminus s \in \varOmega^{\rm{S+}}, \forall r, \forall h, \label{GandTEP26}\\
&0\leq p_{srh}^{\rm{S^{C}}}\leq m_{s}^{\rm{S}}\myol{P}_{s}^{\rm{S}}, \quad \forall s \in\varOmega^{\rm{S+}}, \forall r, \forall h, \label{GandTEP27}\\
&0\leq p_{srh}^{\rm{S^{D}}}\leq \myol{P}_{s}^{\rm{S}}, \quad \forall s \setminus s \in \varOmega^{\rm{S+}}, \forall r, \forall h, \label{GandTEP28}\\
&0\leq p_{srh}^{\rm{S^{D}}}\leq m_{s}^{\rm{S}}\myol{P}_{s}^{\rm{S}}, \quad \forall s \in\varOmega^{\rm{S+}}, \forall r, \forall h, \label{GandTEP29}\\
&0\leq p_{wrh}^{\rm{W}}\leq \alpha_{wrh}\myol{P}_{w}^{\rm{W}}, \quad \forall w \setminus w \in \varOmega^{\rm{W+}}, \forall r, \forall h, \label{GandTEP30}\\
&0\leq p_{wrh}^{\rm{W}}\leq \alpha_{wrh}\overline{p}_{w}^{\rm{W}}, \quad \forall w\in\varOmega^{\rm{W+}}, \forall r, \forall h, \label{GandTEP31}\\
&\theta_{nrh} = 0, \quad n:ref., \quad \forall r, \forall h, \label{GandTEP32}
\end{align}
\end{subequations}
where variables in set $\Phi = \left\{ m_{s}^{\rm{S}}\right.$, $\forall s \in \varOmega^{\rm{S+}}$; $p^{\rm{G}}_{grh}$, $\forall g$, $\forall r$, $\forall h$; $\overline{p}^{\rm{G}}_{g}$, $\forall g \in \varOmega^{\rm{G+}}$; $p^{\rm{L}}_{\ell rh}$, $\forall \ell$, $\forall r$, $\forall h$; $p^{\rm{LS}}_{drh}$, $\forall d$, $\forall r$, $\forall h$; $e_{srh}^{\rm{S}}$, $p^{\rm{S^{C}}}_{srh}$, $p^{\rm{S^{D}}}_{srh}$, $\forall s$, $\forall r$, $\forall h$; $p^{\rm{W}}_{wrh}$, $\forall w$, $\forall r$, $\forall h$; $\overline{p}^{\rm{W}}_{w}$, $\forall w \in \varOmega^{\rm{W+}}$; $x^{\rm{L}}_{\ell}$, $\forall \ell \in \varOmega^{\rm{L+}}$; $\theta_{nrh}$, $\forall n$, $\forall r$, $\forall h \left. \right\}$ are the optimization variables of problem \eqref{GandTEP}.

The objective function \eqref{GandTEP1} represents the aim of the G\&TEP problem, which is minimizing the expansion (generation, storage, and transmission facilities) and operation (power produced by conventional generating units and load-shedding) costs.
The terms associated with operation costs are multiplied by the weight of the corresponding representative day, $\sigma_{r}$, to make them comparable with expansion costs.
Note that the sum of $\sigma_{r}$ for all the representative days is equal to 365, i.e., the total number of days in a year.

Constraints \eqref{GandTEP2} limit the number of units to be built of each candidate storage facility.
Constraints \eqref{GandTEP3} define $m_{s}^{\rm{S}}$, $\forall s$, as integer variables.
Constraints \eqref{GandTEP4}-\eqref{GandTEP5} impose bounds on the capacity to be built of conventional and wind-power generating units, respectively.
Constraints \eqref{GandTEP6} define $x_{\ell}^{\rm{L}}$ as binary variables that indicate whether a candidate transmission line is built ($x_{\ell}^{\rm{L}}=1$) or not ($x_{\ell}^{\rm{L}}=0$).
Constraints \eqref{GandTEP7}-\eqref{GandTEP10} impose investment budgets for building candidate conventional generating units, transmission lines, storage, and wind-power units, respectively.
Constraints \eqref{GandTEP11}-\eqref{GandTEP31} are the operation constraints and comprise equations \eqref{GandTEP11} that impose the generation-demand balance at each node, where demand factors $\beta_{drh}$, $\forall d$, $\forall r$, $\forall h$, are linked to the output of the K-means method described in Section \ref{SecMeth}; constraints \eqref{GandTEP12}-\eqref{GandTEP13} that define the power flows through existing and candidate transmission lines, respectively, which are limited by constraints \eqref{GandTEP14}-\eqref{GandTEP15}; equations \eqref{GandTEP16} that define the energy stored in storage units for all representative days and hours excluding the first hour of each day; equations \eqref{GandTEP17}-\eqref{GandTEP18} that define the energy stored in existing and candidate storage units, respectively, for the first hour of all representative days; constraints \eqref{GandTEP19}-\eqref{GandTEP20} which ensure that existing and candidate storage units, respectively, store a minimum amount of energy at the end of each representative day; constraints \eqref{GandTEP21}-\eqref{GandTEP22} that impose bounds on the energy stored in the existing and candidate storage units, respectively; constraints \eqref{GandTEP23}-\eqref{GandTEP24} that impose bounds on the power produced by existing and candidate conventional generating units, respectively; constraints \eqref{GandTEP25} that limit the load shed of demands; constraints \eqref{GandTEP26}-\eqref{GandTEP27} that impose bounds on the charging power of existing and candidate storage units, respectively; constraints \eqref{GandTEP28}-\eqref{GandTEP29} that  impose bounds on the discharging power of existing and candidate storage units, respectively; constraints \eqref{GandTEP30}-\eqref{GandTEP31} that impose bounds on the power produced by existing and candidate wind-power units, respectively, where wind-power capacity factors $\alpha_{wrh}$, $\forall w$, $\forall r$, $\forall h$, are associated with the output of the K-means method described in Section \ref{SecMeth}; and constraints \eqref{GandTEP32} which define the voltage angle at the reference node.

It is important to mention that the network constraints are modeled in the G\&TEP problem using a DC model without losses for the sake of simplicity. 
In addition, fixed costs are not considered and the capacity to be installed of each generating unit, i.e., variables $\overline{p}_{g}^{\rm{G}}$, $\forall g \in \varOmega^{\rm{G+}}$, are considered continuous.

The G\&TEP problem \eqref{GandTEP} is a mixed-integer nonlinear programming (MINLP) model.
Nonlinear terms are $x_{\ell}^{{\rm{L}}} \theta_{n r h }$ in constraints \eqref{GandTEP13}, i.e., products of binary and continuous variables.
These nonlinear terms can be replaced by exact equivalent mixed-integer linear expressions as explained, e.g., in \cite{Floudas95}.
Thus, the G\&TEP problem \eqref{GandTEP} can be finally formulated as a mixed-integer programming (MILP) model that can be solved using available branch-and-cut solvers, e.g., CPLEX \cite{Cplex2016}.

\section{Case Study} \label{SecCS}

\subsection{Data}

\begin{figure}[H]
	\centering
	\includegraphics[width=13 cm]{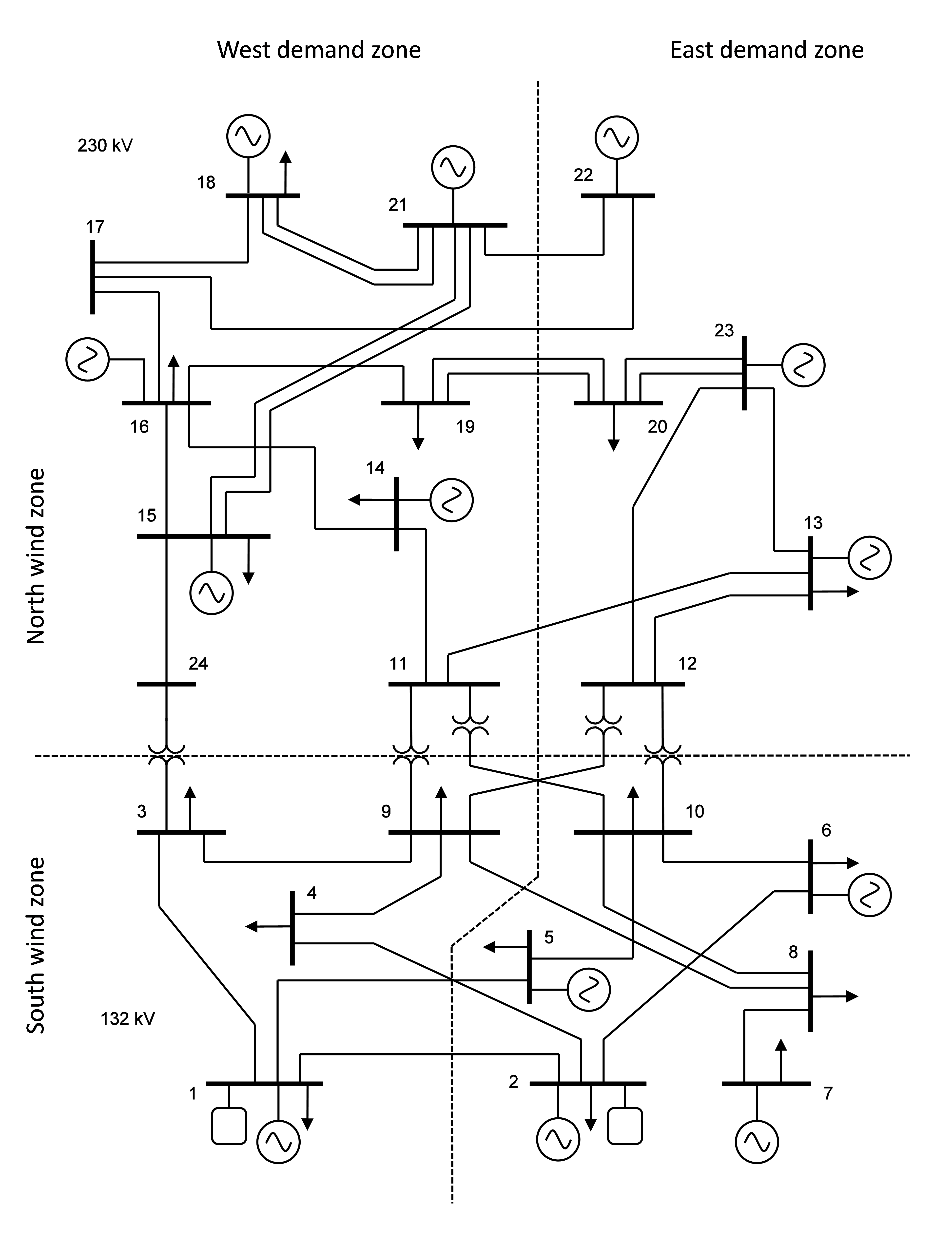}
	\caption{Modified version of the IEEE RTS.}
	\label{f501}
\end{figure}

We apply the expansion model described in Section \ref{SecForm} to the modified version of the IEEE Reliability Test System (RTS) \cite{RTS} that is depicted in Fig. \ref{f501}.
This electric energy system comprises 11 conventional generating units, 17 demands, 24 nodes, two storage units, 38 transmission lines and two wind-power units.
Table \ref{t501} provides the conventional generating unit data; Table \ref{t502} supplies the demand data; storage unit data is presented in Table \ref{t503}; the transmission line data can be consulted in Table \ref{t504}; and Table \ref{t505} provides the wind-power unit data.
It is necessary to mention that the annualized investment costs of candidate storage units, which are showed in Table  \ref{t503}, are based on the data collected in \cite{FernandezBlanco2017}.
We consider a set of values taking the average value of the costs provided in the two scenarios considered in \cite{FernandezBlanco2017}, as it is displayed in equation \eqref{IS}.
\begin{equation}
I_{s}^{\rm{S}}=60,000\myol{E}_{s}^{\rm{S}}+1,000,000\myol{P}_{s}^{\rm{S}} \label{IS}
\end{equation}

\begin{figure}[H]
	\centering
	\includegraphics[width=11 cm]{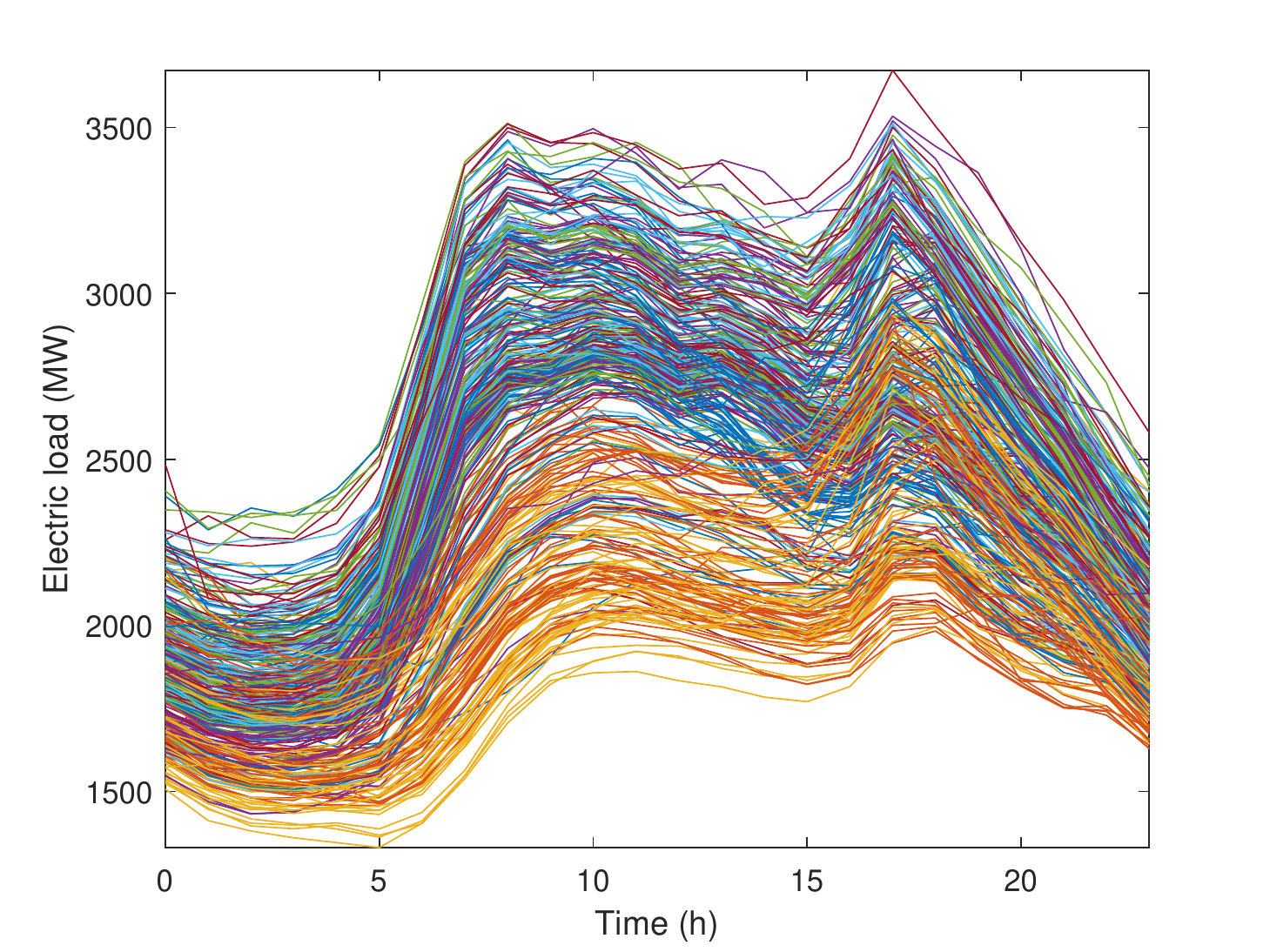}
	\caption{Case study: diary evolution of electric load in the west zone during a year.}
	\label{f502}
\end{figure}

\begin{figure}[H]
	\centering
	\includegraphics[width=11 cm]{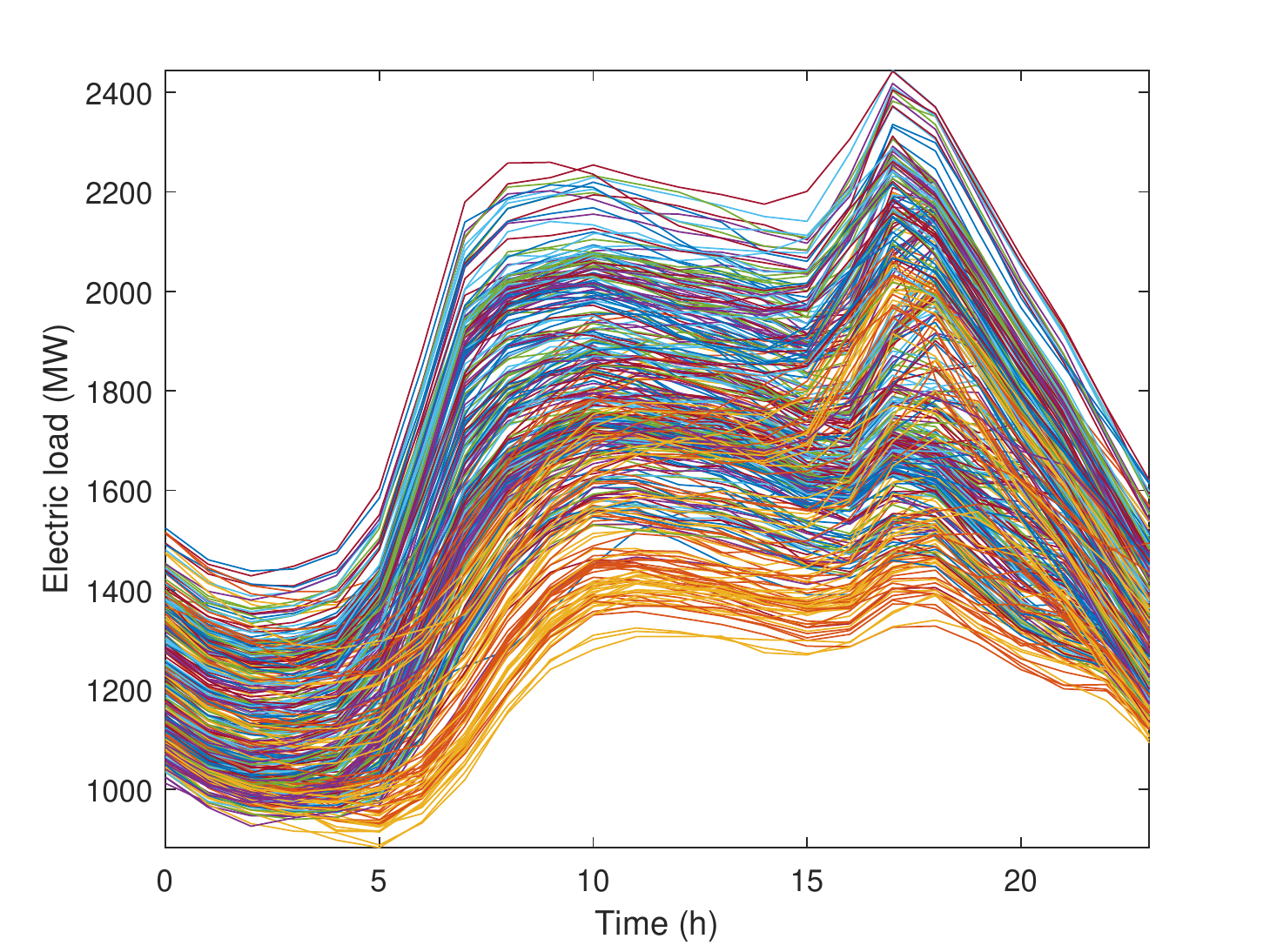}
	\caption{Case study: diary evolution of electric load in the east zone during a year.}
	\label{f503}
\end{figure}

\begin{figure}[H]
	\centering
	\includegraphics[width=10.5 cm]{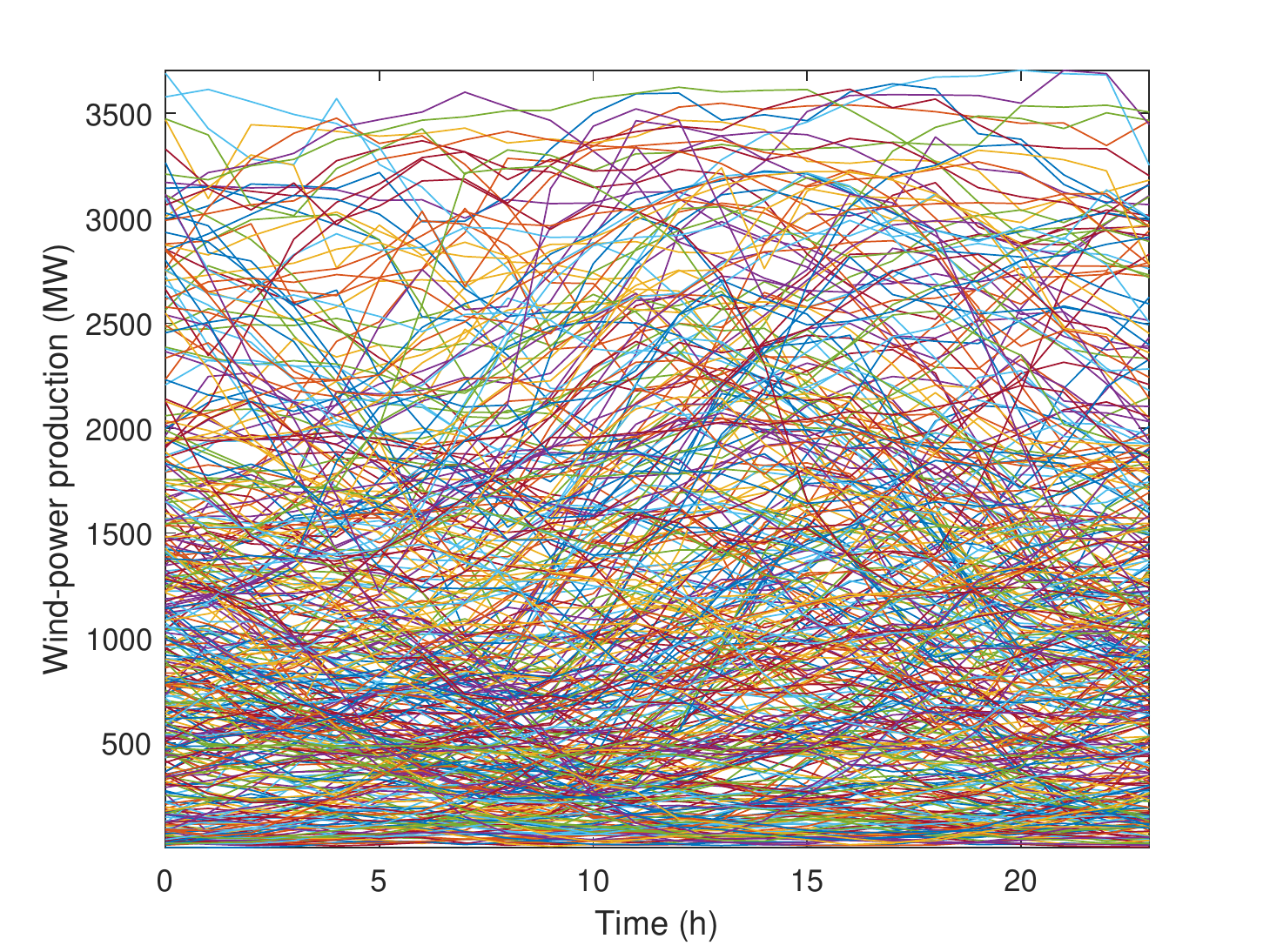}
	\caption{Case study: diary evolution of wind-power production in the north zone during a year.}
	\label{f504}
\end{figure}

\begin{figure}[H]
	\centering
	\includegraphics[width=10.5 cm]{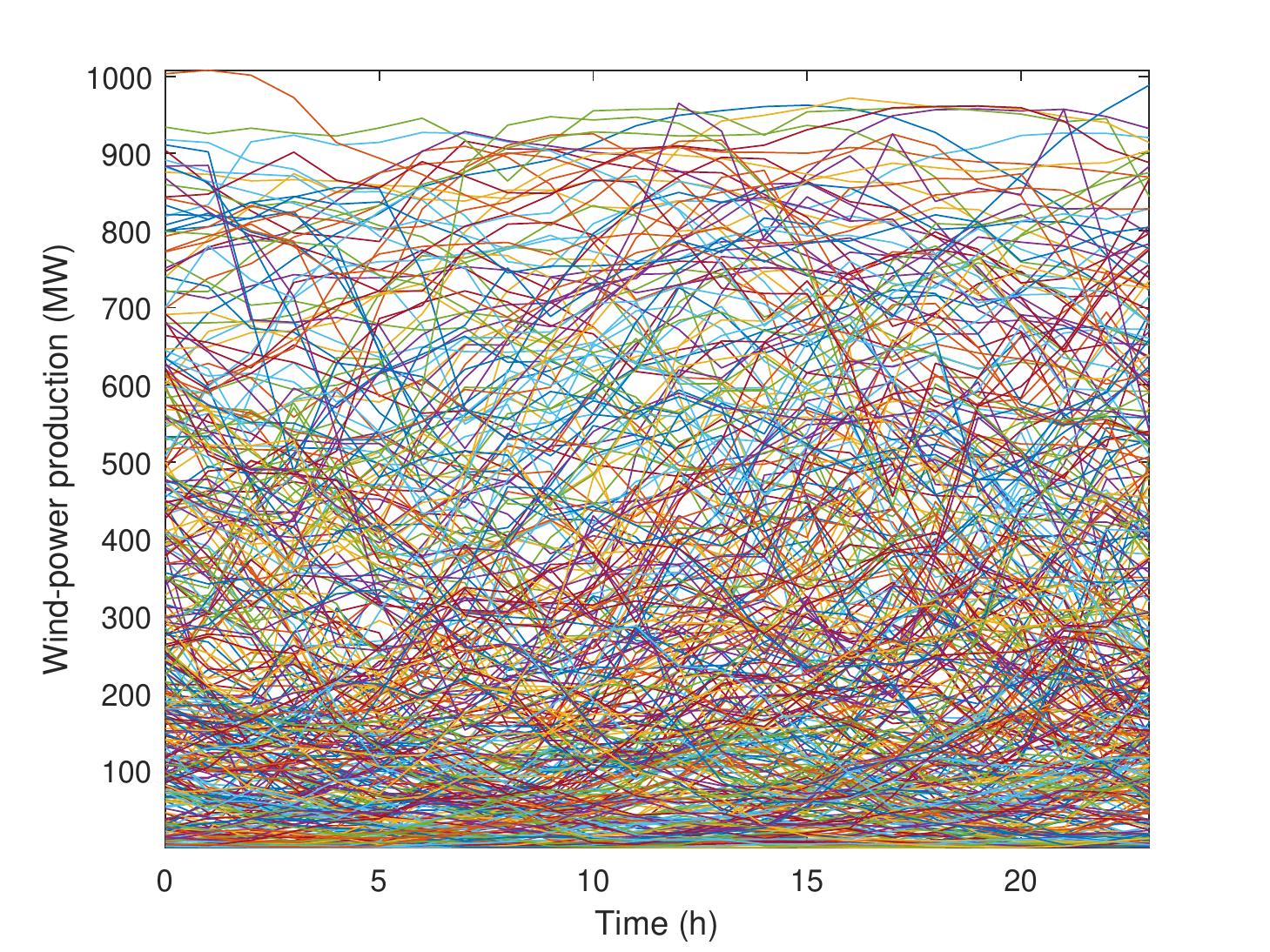}
	\caption{Case study: diary evolution of wind-power production in the south zone during a year.}
	\label{f505}
\end{figure}

We consider that wind-power production and electric load can change their values depending on the zone of the electric energy system where wind-power units and demands are located.
On the one hand, demands are allocated to the west and east zones of the system, as illustrated in Figs. \ref{f502} and \ref{f503}.
On the other hand, wind-power units are distributed between the north and south zones of the system, as depicted in Figs. \ref{f504} and \ref{f505}.
As in Section \ref{SecMeth}, the historical data of electric load and wind-power production have been acquired from \cite{Energinet}.
It is remarkable to mention that the peak values of electric load in the west zone are greater than in the east zone.
In addition, the maximum values of wind-power production are associated with the north zone.
It is expected that the need to supply the high demands in the west zone will condition the investment decision making of the expansion problem.

It is supposed that we work with hourly data, thus the duration of time steps, $\Delta \tau$, is equal to one hour.
We consider that the charging and discharging efficiency of storage units is equal to 90 \%.
The energy initially stored in storage units is assumed to be zero for all the representative days.
Node 1 is the reference node of the optimization problem.
The parameter $F$ receives a value of 500,000.
Due to the presence of transformers in the electric energy system considered as it can be noticed in Fig. \ref{f501}, we select a base power of 100 MW.
It is supposed that the values of the parameters $\alpha_{wrh}$ and $\beta_{drh}$ for each representative day and hour are the same for all the wind-power units and demands, respectively. Both parameters are obtained from the K-means methods.

Instead of considering a different investment budget for the building of each candidate generating/storage unit or transmission line, we consider a total investment budget, $\myol{I}^{\rm{T}}$, which is distributed among the different types of facilities.
Thus, it is supposed that constraints \eqref{GandTEP7}-\eqref{GandTEP10} of problem \eqref{GandTEP} are replaced by constraint \eqref{IT} from now on.
Therefore, we consider a total investment budget of \$2,000 million.
The annualized investment costs are 10 \% of the total costs.
\begin{equation}
\sum_{g\in\varOmega^{\rm{G+}}}I_{g}^{\rm{G}}\overline{p}_{g}^{\rm{G}}+\sum_{\ell\in\varOmega^{\rm{L+}}}I_{\ell}^{\rm{L}}x_{\ell}^{\rm{L}}+\sum_{s\in\varOmega^{\rm{S+}}}I_{s}^{\rm{S}}m_{s}^{\rm{S}}+\sum_{w\in\varOmega^{\rm{W+}}}I_{w}^{\rm{W}}\overline{p}_{w}^{\rm{W}}\leq \myol{I}^{\rm{T}} \label{IT}
\end{equation}
\begin{longtable}{|c|c|c|c|c|c|} \hline
	Conventional generating unit&Node&$\myol{P}_{g}^{\rm{G}}$ [MW]&$C_{g}^{\rm{G}}$ [\$/MWh]&$\tilde{I}_{g}^{\rm{G}}$ [\$/MW] \\ \hline 
	$g_{1}$&$n_{1}$&172&75&-\\ \hline
	$g_{2}$&$n_{2}$&172&77&-\\ \hline
	$g_{3}$&$n_{7}$&240&75&-\\ \hline
	$g_{4}$&$n_{13}$&285&70&-\\ \hline
	$g_{5}$&$n_{14}$&200&72&-\\ \hline
	$g_{6}$&$n_{15}$&215&67&-\\ \hline
	$g_{7}$&$n_{16}$&155&69&-\\ \hline
	$g_{8}$&$n_{18}$&400&71&-\\ \hline
	$g_{9}$&$n_{21}$&400&68&-\\ \hline
	$g_{10}$&$n_{22}$&300&70&-\\ \hline
	$g_{11}$&$n_{23}$&260&65&-\\ \hline
	$g_{12}$&$n_{3}$&250&55&100,000\\ \hline
	$g_{13}$&$n_{8}$&250&53&100,000\\ \hline
	$g_{14}$&$n_{9}$&200&60&100,000\\ \hline
	$g_{15}$&$n_{12}$&200&58&100,000\\ \hline
	$g_{16}$&$n_{16}$&250&54&100,000\\ \hline
	$g_{17}$&$n_{19}$&200&59&100,000\\ \hline
	$g_{18}$&$n_{20}$&250&55&100,000\\\hline            
	\multicolumn{5}{c}{}\\
	\caption{Case study: conventional generating unit data.}
	\label{t501}
\end{longtable}
\begin{longtable}{|c|c|c|c|c|c|} \hline
	Demand&Node&Zone&$\myol{P}_{d}^{\rm{D}}$ [MW]&$C_{d}^{\rm{LS}}$ [\$/MWh] \\ \hline 
	$d_{1}$&$n_{1}$&West&270.0&30,000\\ \hline
	$d_{2}$&$n_{2}$&East&242.5&30,000\\ \hline
	$d_{3}$&$n_{3}$&West&450.0&30,000\\ \hline
	$d_{4}$&$n_{4}$&West&185.0&30,000\\ \hline
	$d_{5}$&$n_{5}$&East&177.5&30,000\\ \hline
	$d_{6}$&$n_{6}$&East&340.0&30,000\\ \hline
	$d_{7}$&$n_{7}$&East&312.5&30,000\\ \hline
	$d_{8}$&$n_{8}$&East&427.5&30,000\\ \hline
	$d_{9}$&$n_{9}$&West&437.5&30,000\\ \hline
	$d_{10}$&$n_{10}$&East&487.5&30,000\\ \hline
	$d_{11}$&$n_{13}$&East&662.5&30,000\\ \hline
	$d_{12}$&$n_{14}$&West&485.0&30,000\\ \hline
	$d_{13}$&$n_{15}$&West&792.5&30,000\\ \hline
	$d_{14}$&$n_{16}$&West&250.0&30,000\\ \hline
	$d_{15}$&$n_{18}$&West&832.5&30,000\\ \hline
	$d_{16}$&$n_{19}$&West&452.5&30,000\\ \hline
	$d_{17}$&$n_{20}$&East&320.0&30,000\\\hline         
	\multicolumn{5}{c}{}\\    
	\caption{Case study: demand data.}
	\label{t502}
\end{longtable}
\begin{longtable}{|c|c|c|c|c|c|c|} \hline
	Storage unit&Node&$\myol{M}_{s}^{\rm{S}}$&$\myol{E}_{s}^{\rm{S}}$ [MWh]&$\myol{P}_{s}^{\rm{S}}$ [MW]&$\tilde{I}_{s}^{\rm{S}}$ [\$] \\ \hline 
	$s_{1}$&$n_{1}$&-&100&50&-\\ \hline
	$s_{2}$&$n_{2}$&-&100&50&-\\ \hline
	$s_{3}$&$n_{13}$&2&250&125&14,000,000\\ \hline
	$s_{4}$&$n_{15}$&3&250&125&14,000,000\\ \hline
	$s_{5}$&$n_{18}$&2&200&100&11,200,000\\ \hline
	$s_{6}$&$n_{21}$&1&300&150&16,800,000\\ \hline
	$s_{7}$&$n_{23}$&1&400&200&22,400,000\\ \hline
	\multicolumn{7}{c}{}\\
	\caption{Case study: storage unit data.}	
	\label{t503}
\end{longtable}
\newpage
\begin{longtable}{|c|c|c|c|c|c|c|} \hline 
	Transmission line&From bus&To bus&$1 / B_{\ell}$ [pu]&$\myol{P}_{\ell}^{\rm{L}}$ [MW]&$\tilde{I}_{\ell}^{\rm{L}}$ [\$] \\ \hline 
	$\ell_{1}$&$n_{1}$&$n_{2}$&0.014&150&-\\ \hline
	$\ell_{2}$&$n_{1}$&$n_{3}$&0.211&150&-\\ \hline
	$\ell_{3}$&$n_{1}$&$n_{5}$&0.085&150&-\\ \hline
	$\ell_{4}$&$n_{2}$&$n_{4}$&0.127&150&-\\ \hline
	$\ell_{5}$&$n_{2}$&$n_{6}$&0.192&150&-\\ \hline
	$\ell_{6}$&$n_{3}$&$n_{9}$&0.119&150&-\\ \hline
	$\ell_{7}$&$n_{3}$&$n_{24}$&0.084&150&-\\ \hline
	$\ell_{8}$&$n_{4}$&$n_{9}$&0.104&150&-\\ \hline
	$\ell_{9}$&$n_{5}$&$n_{10}$&0.088&150&-\\ \hline
	$\ell_{10}$&$n_{6}$&$n_{10}$&0.061&150&-\\ \hline
	$\ell_{11}$&$n_{7}$&$n_{8}$&0.061&150&-\\ \hline
	$\ell_{12}$&$n_{8}$&$n_{9}$&0.161&150&-\\ \hline
	$\ell_{13}$&$n_{8}$&$n_{10}$&0.165&150&-\\ \hline
	$\ell_{14}$&$n_{9}$&$n_{11}$&0.084&150&-\\ \hline
	$\ell_{15}$&$n_{9}$&$n_{12}$&0.084&150&-\\ \hline
	$\ell_{16}$&$n_{10}$&$n_{11}$&0.084&150&-\\ \hline
	$\ell_{17}$&$n_{10}$&$n_{12}$&0.084&150&-\\ \hline
	$\ell_{18}$&$n_{11}$&$n_{13}$&0.048&150&-\\ \hline
	$\ell_{19}$&$n_{12}$&$n_{14}$&0.042&150&-\\ \hline
	$\ell_{20}$&$n_{12}$&$n_{13}$&0.048&150&-\\ \hline
	$\ell_{21}$&$n_{12}$&$n_{23}$&0.087&150&-\\ \hline
	$\ell_{22}$&$n_{13}$&$n_{23}$&0.075&150&-\\ \hline
	$\ell_{23}$&$n_{14}$&$n_{16}$&0.059&150&-\\ \hline
	$\ell_{24}$&$n_{15}$&$n_{16}$&0.017&150&-\\ \hline
	$\ell_{25}$&$n_{15}$&$n_{21}$&0.049&150&-\\ \hline
	$\ell_{26}$&$n_{15}$&$n_{21}$&0.049&150&-\\ \hline
	$\ell_{27}$&$n_{15}$&$n_{24}$&0.052&150&-\\ \hline
	$\ell_{28}$&$n_{16}$&$n_{17}$&0.026&150&-\\ \hline
	$\ell_{29}$&$n_{16}$&$n_{19}$&0.023&150&-\\ \hline
	$\ell_{30}$&$n_{17}$&$n_{18}$&0.014&150&-\\ \hline
	$\ell_{31}$&$n_{17}$&$n_{22}$&0.105&150&-\\ \hline
	$\ell_{32}$&$n_{18}$&$n_{21}$&0.026&150&-\\ \hline
	$\ell_{33}$&$n_{18}$&$n_{21}$&0.026&150&-\\ \hline
	$\ell_{34}$&$n_{19}$&$n_{20}$&0.040&150&-\\ \hline
	$\ell_{35}$&$n_{19}$&$n_{20}$&0.040&150&-\\ \hline
	$\ell_{36}$&$n_{20}$&$n_{23}$&0.220&150&-\\ \hline
	$\ell_{37}$&$n_{20}$&$n_{23}$&0.220&150&-\\ \hline
	$\ell_{38}$&$n_{21}$&$n_{22}$&0.068&150&-\\ \hline
	$\ell_{39}$&$n_{2}$&$n_{7}$&0.120&175&106,670\\ \hline
	$\ell_{40}$&$n_{6}$&$n_{13}$&0.140&175&113,330\\ \hline
	$\ell_{41}$&$n_{7}$&$n_{8}$&0.165&175&111,000\\ \hline
	$\ell_{42}$&$n_{11}$&$n_{19}$&0.048&500&228,940\\ \hline
	$\ell_{43}$&$n_{11}$&$n_{24}$&0.048&500&228,940\\ \hline
	$\ell_{44}$&$n_{12}$&$n_{19}$&0.075&500&416,250\\ \hline
	\multicolumn{7}{c}{}\\
	\caption{Case study: transmission line data.}
	\label{t504}
\end{longtable}
\begin{longtable}{|c|c|c|c|c|c|} \hline
	Wind-power unit&Node&Zone&$\myol{P}_{w}^{\rm{W}}$ [MW]&$\tilde{I}_{w}^{\rm{W}}$ [\$/MW] \\ \hline 
	$r_{1}$&$n_{5}$&South&200&-\\ \hline
	$r_{2}$&$n_{6}$&South&200&-\\ \hline
	$r_{3}$&$n_{7}$&South&300&300,000\\ \hline
	$r_{4}$&$n_{10}$&South&400&300,000\\ \hline
	$r_{5}$&$n_{16}$&North&300&300,000\\ \hline
	$r_{6}$&$n_{20}$&North&300&300,000\\ \hline
	\multicolumn{5}{c}{}\\
	\caption{Case study: wind-power unit data.}
	\label{t505}
\end{longtable}

\subsection{Results}

First of all, we solve the G\&TEP problem using all the historical data to find the exact solution in order to compare it with the results obtained using representative days provided by both K-mean methods.
However, it is necessary to make some changes in the formulation of problem \eqref{GandTEP} to properly characterize the continuity in time of the historical data.
Thus, constraints \eqref{GandTEP17}-\eqref{GandTEP18} are replaced by constraints \eqref{ES1}, which allude to the energy stored in each storage unit during the first hour of all the days except the first one relating it to the energy stored in the same storage unit during the last hour of the previous day.
In addition, constraints \eqref{GandTEP19}-\eqref{GandTEP20} are replaced by constraints \eqref{ES2}-\eqref{ES3}, which refer to the energy stored in each existing and candidate storage unit, respectively, during the first hour of the first day linking it to the energy initially stored in the same storage unit in the first day, $E_{sr_{0}}^{\rm{S}}$.
\begin{align}
&e_{srh_{1}}^{\rm{S}}=e_{sr-1h_{24}}^{\rm{S}}+\left(  p_{srh_{1}}^{\rm{S^{C}}}\eta_{s}^{\rm{S^{C}}}-\frac{p_{srh_{1}}^{\rm{S^{D}}}}{\eta_{s}^{\rm{S^{D}}}}\right) \varDelta\tau \quad \forall s, \forall r > 1 \label{ES1}\\
&e_{sr_{1}h_{1}}^{\rm{S}}=E_{sr_{0}}^{\rm{S}}+\left(  p_{sr_{1}h_{1}}^{\rm{S^{C}}}\eta_{s}^{\rm{S^{C}}}-\frac{p_{sr_{1}h_{1}}^{\rm{S^{D}}}}{\eta_{s}^{S^{D}}}\right) \varDelta\tau \quad \forall s \setminus s \in \varOmega^{\rm{S+}} \label{ES2}\\
&e_{sr_{1}h_{1}}^{\rm{S}}=m_{s}^{\rm{S}}E_{sr_{0}}^{\rm{S}}+\left(  p_{sr_{1}h_{1}}^{\rm{S^{C}}}\eta_{s}^{\rm{S^{C}}}-\frac{p_{sr_{1}h_{1}}^{\rm{S^{D}}}}{\eta_{s}^{S^{D}}}\right) \varDelta\tau \quad \forall s \in \varOmega^{\rm{S+}} \label{ES3}
\end{align}

Having made these changes, the G\&TEP problem is solved using the 366 days of historical data, due to the fact that the year considered is a leap year. The total annual cost obtained, $CT$, amounts to \$3,124 million. The results show that the 0.14 \% of the total demand is not supplied. The computation time required to obtain the exact solution is 55 h 28 min.

The steps that should be followed in order to make the results obtained using representative days comparable with the exact solution are presented below:
\begin{itemize}
	\item Step 1: Solve the G\&TEP problem using representative days obtained applying the clustering methods.
	\item Step 2: Fix the values of the decision variables ($m_{s}^{\rm{S}}$, $\forall s \in \varOmega^{\rm{S+}}$; $\overline{p}_{g}^{\rm{G}}$, $\forall g \in \varOmega^{\rm{G+}}$; $\overline{p}_{w}^{\rm{W}}$, $\forall w \in \varOmega^{\rm{W+}}$; $x_{\ell}^{\rm{L}}$, $\forall \ell \in \varOmega^{\rm{L+}}$) obtained in Step 1 and solve the G\&TEP problem using all the historical data.
	\item Step 3: Calculate the percent error, $\varepsilon_{CT}$, associated with the total annual cost obtained in Step 2, $CT^{\rm{K}}$, with regard to the  total annual cost provided by the exact solution, $CT^{\rm{E}}$, applying the equation \eqref{error}.
	\begin{equation}
	\varepsilon_{CT} = \frac{\lvert CT^{\rm{K}}-CT^{\rm{E}}\rvert}{CT^{\rm{E}}}\cdot 100 \label{error}
	\end{equation}
\end{itemize}

These steps are followed in the case study using a set of values of the parameter $K$ ranging from 10 to 80, being 366 the maximum value which could be selected.
It means that we work with an equivalent amount of data ranging from 3 to 22 \% of all the historical data considered.

\begin{figure}[H]
	\centering
	\includegraphics[width=8 cm]{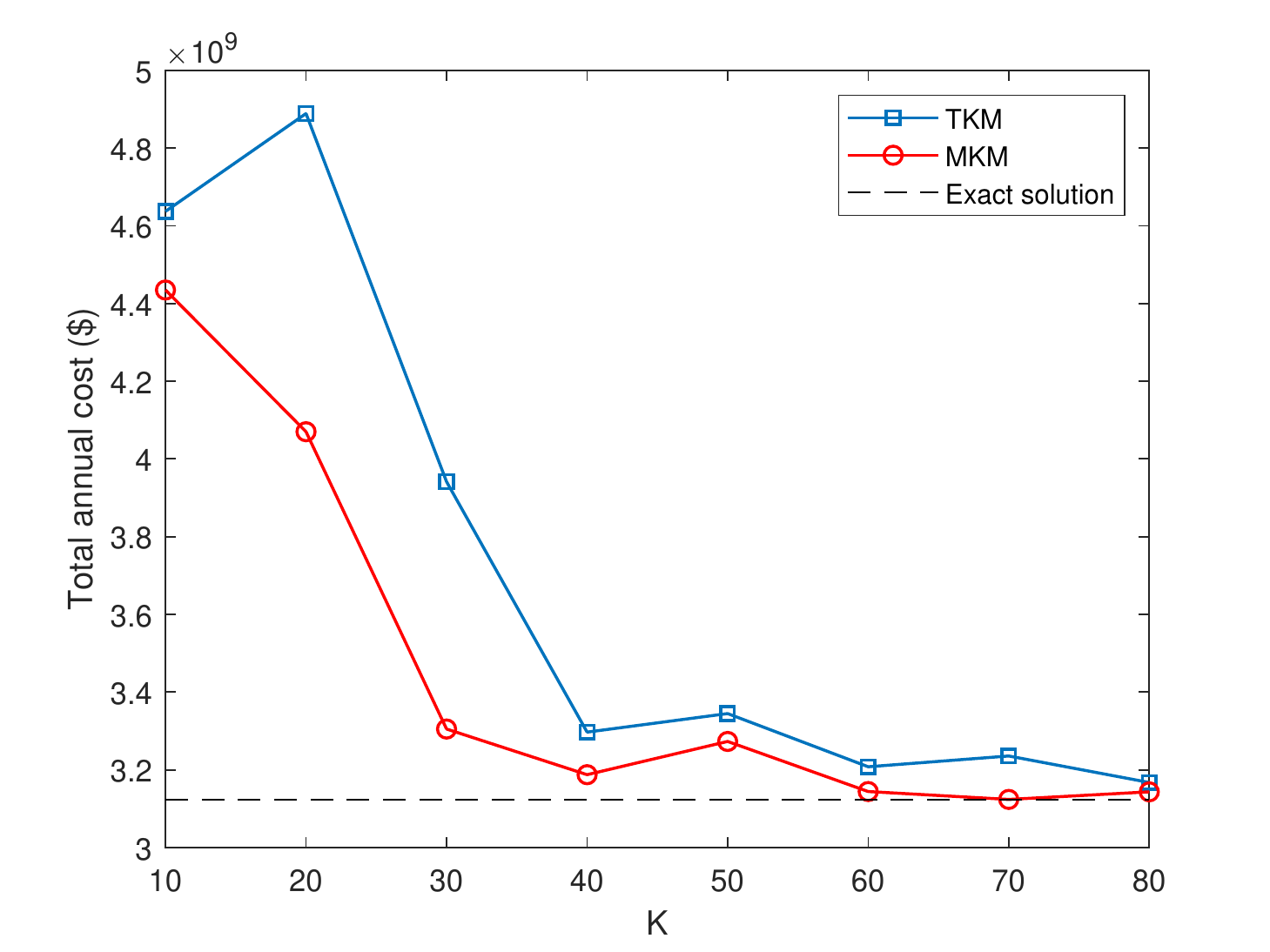}
	\caption{Case study: total annual cost obtained using different values of $K$ and clustering methods.}
	\label{f506}
\end{figure}

Fig. \ref{f506} depicts the total annual cost obtained using different values of $K$ and clustering methods.
We observe that the MKM presents values of the total annual cost closer to the exact solution than those obtained using the TKM for all the cases evaluated.
Note that the differences among the results obtained using the clustering methods and the exact solution generally decrease at the same time that the value of $K$ increases.
However, this is not always true because, for instance, this differences are greater considering $K=50$ than in the case of using $K=40$.
Due to the high total investment budget taken into account, most of the candidate facilities considered are built and more than the 99 \% of the demand is supplied.

\begin{figure}[H]
	\centering
	\includegraphics[width=8 cm]{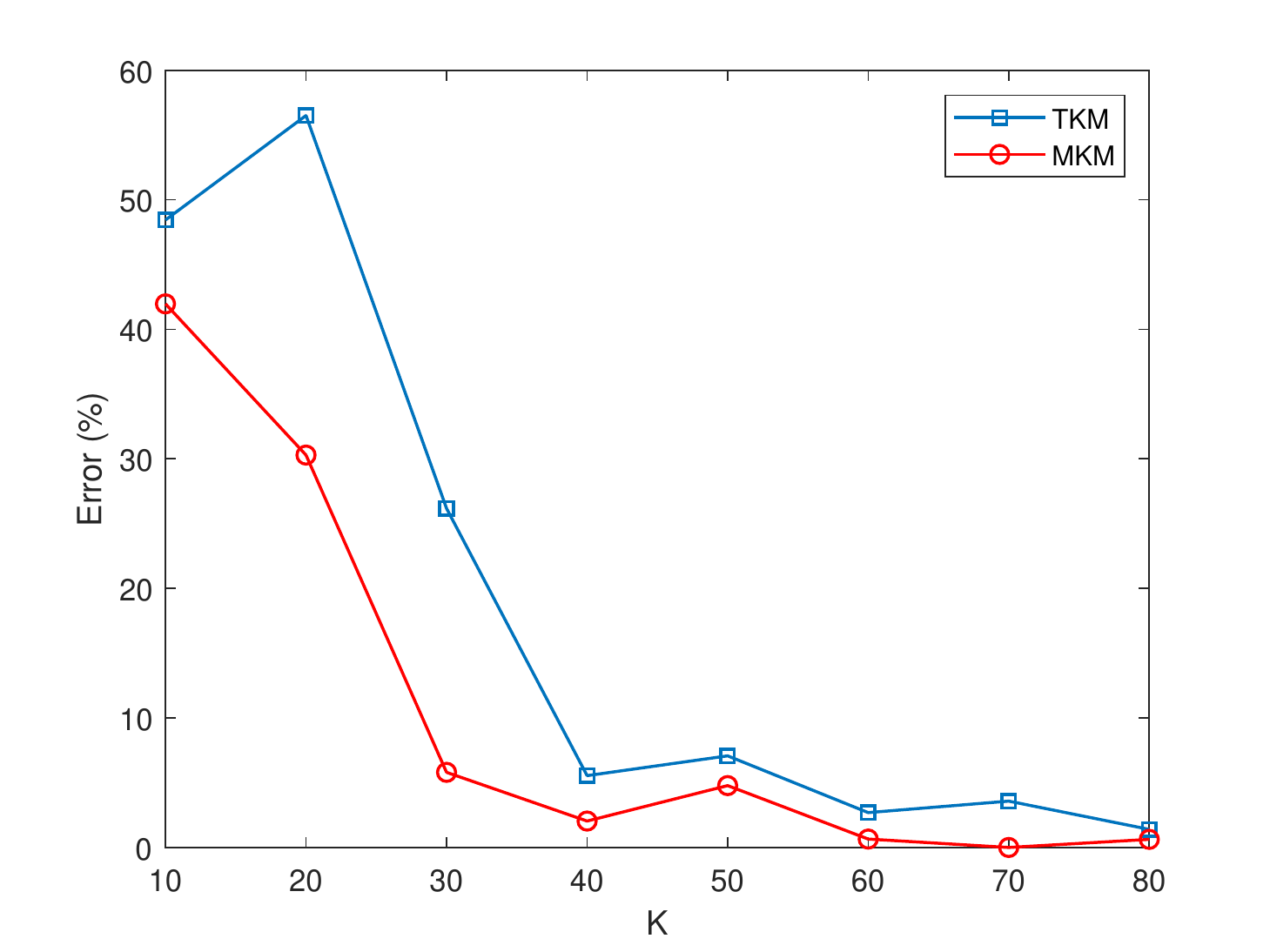}
	\caption{Case study: error of the total annual cost obtained using different values of $K$ and clustering methods.}
	\label{f507}
\end{figure}

Fig. \ref{f507} illustrates the error of the total annual cost obtained using different values of $K$ and clustering methods.
It is clear that the MKM provides results with less error than those obtained using the TKM for all the cases analyzed, especially in those where the parameter $K$ presents a low value.

\begin{figure}[H]
	\centering
	\includegraphics[width=8 cm]{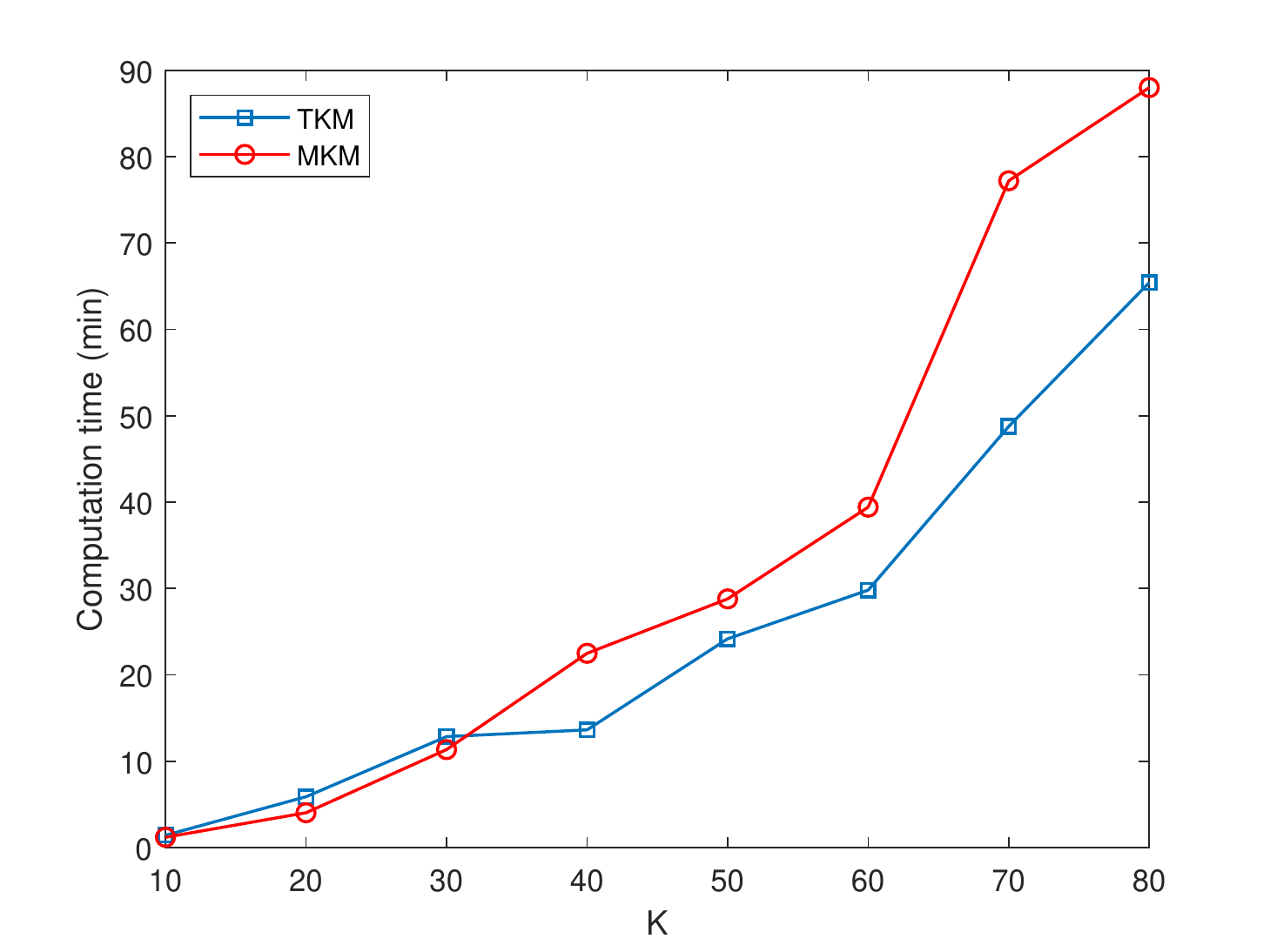}
	\caption{Case study: computation times obtained using different values of $K$ and clustering methods.}
	\label{f508}
\end{figure}

Although it is fundamental to determine which clustering method provides the closest results to the exact solution, we should also analyze the computation times, obtained in Step 1 of the process described above, in the cases under study.
It is relevant in Fig. \ref{f508} that the TKM generally provides shorter computation times, especially in those cases where the parameter $K$ presents a high value.
However, it should be taken into account that the possible saturation of the server used to solve the G\&TEP problem, caused by its concurrent use, may have influenced in the values of the computation times obtained.
In addition, note that there is a rising trend of the computation times as well as it is increased the value of $K$.
The result of Figs. \ref{f506}, \ref{f507} and \ref{f508} are collected in Table \ref{t506}.

\begin{table}[H]
	\centering
	\begin{tabular}{|c|c|c|c|c|c|c|}
		\cline{2-7}
		\multicolumn{1}{c|}{}&\multicolumn{2}{c|}{\multirow{2}{*}{CT [$\cdot 10^{9} \; \$ $]}}&\multicolumn{2}{|c|}{\multirow{2}{*}{$\varepsilon_{CT} (\%)$}}&\multicolumn{2}{c|}{Computation}\\
		\multicolumn{1}{c|}{}&\multicolumn{1}{c}{}&\multicolumn{1}{c|}{}&\multicolumn{1}{c}{}&\multicolumn{1}{c|}{}&\multicolumn{2}{c|}{time [min]}\\ \hline
		$K$&TKM&MKM&TKM&MKM&TKM&MKM\\ \hline
		10&4.69&4.43&48.44&41.97&1&1\\ \hline
		20&4.89&4.07&56.51&30.29&6&4\\ \hline
		30&3.94&3.30&26.18&5.80&13&11\\ \hline
		40&3.30&3.19&5.55&2.04&14&22\\ \hline
		50&3.34&3.27&7.08&4.78&24&29\\ \hline
		60&3.21&3.14&2.70&0.66&30&39\\ \hline
		70&3.24&3.12&3.58&0.01&49&77\\ \hline
		80&3.17&3.14&1.40&0.63&65&88\\ \hline
	\end{tabular}
	\caption{Case study: analysis of the result obtained from the clustering methods.}
	\label{t506}
\end{table}

Taking into account the results commented before, we consider that the MKM provides better results than the TKM, especially regarding  the error of the total annual cost. Although the computation times obtained using the MKM are generally greater than those acquired using the TKM, in several of the cases evaluated the error provided by the MKM in a given time is less than the error obtained using the TKM and the same amount of time.
For instance, the MKM presents a 2.04 \% of error using 40 representative days in 22 min, while the TKM spends 30 min to obtain a 2.70 \% of error using 60 representative days.
Due to this and the possible saturation problems in the server mentioned before, we consider that the results associated with the error are more relevant than those linked to the computation times.

\subsection{Computation Times}

The results of this case study are obtained using CPLEX \cite{Cplex2016} under GAMS \cite{GAMS} on an Intel Xeon E7-4820 computer with 4 processors at 2 GHz and 128 GB of RAM.

The computation time required to obtain the exact solution is 55 h 28 min.
Regarding the resolution of the G\&TEP problem using representative days, the corresponding computation times are collected in Table \ref{t506}.

\section{Conclusions} \label{SecConc}

This paper proposes a new clustering method to adequately characterize the maximum and minimum values of the input data.
In addition, we arrange the operating conditions obtained using the K-means method into representative days in order to depict the chronology of the historical data.
This allows us to include storage units in the expansion model considered to solve the G\&TEP problem.

The conclusion of this paper is that the results obtained in the case study using the modified K-means method and different numbers of representative days provide a total annual cost closer to the exact solution than in the case of using the traditional K-means method.
In fact, although the computation times may have been influenced by the saturation of the server used, the results display that in some cases the MKM is able to solve the G\&TEP problem in less time than the TKM using less representative days and achieving a minor error.

%\section*{Acknowledgements}

%The work of L. Baringo has been partially funded by the Ministry of Science of Spain under CICYT Project ENE2015-63879-R (MINECO/FEDER, UE).


\begin{thebibliography}{99}
%
\bibitem{Caramanis82}
M. Caramanis, R. Tabors, K. S. Nochur, and F. Schweppe, ``The introduction of non-dispatchable technologies a decision variables in longterm generation expansion models,'' IEEE Trans. Power App. Syst., vol. PAS-101, no. 8, pp. 2658–2667, 1982.
%
\bibitem{Wogrin13}
S. Wogrin, ``Generation expansion planning in electricity markets with bilevel mathematical programming techniques,'' Ph.D. dissertation, Universidad Pontificia Comillas de Madrid, Madrid, Spain, 2013.
%
\bibitem{Baringo12}
L. Baringo and A. J. Conejo, ``Transmission and wind power investment,'' \emph{IEEE Trans. Power Syst.}, vol. 27, no. 2, pp. 885-893, May 2012.
%
\bibitem{Montoya15}
S. Montoya-Bueno, J. I. Mu\~noz, and J. Contreras, ``A stochastic investment model for renewable generation in distribution systems,'' \emph{IEEE Trans. Sustain. Energy}, vol. 6, no. 4, pp. 1466-1474, Oct. 2015.
%
\bibitem{Weijde12}
A. H. van der Weijde and B. F. Hobbs, ``The economics of planning electricity transmission to accommodate renewables: Using two-stage optimisation to evaluate flexibility and the cost of disregarding uncertainty,'' \emph{Energy Economics}, vol. 34, no. 6, pp. 2089–2101, Nov. 2012.
%
\bibitem{Baringo14}
L. Baringo and A. J. Conejo, ``Strategic wind power investment,'' \emph{IEEE Trans. Power Syst.}, vol. 29, no. 3, pp. 1250-1260, May 2014.
%
\bibitem{Dominguez15}
R. Dom\'inguez, A. J. Conejo, and M. Carri\'on, ``Toward fully renewable electric energy systems,'' \emph{IEEE Trans. Power Syst.}, vol. 30, no. 1, pp. 316-326, Jan. 2015.
%
\bibitem{Baringo13}
L. Baringo and A. J. Conejo, ``Correlated wind-power production and electric load scenarios for investment decisions,'' \emph{Appl. Energy}, vol. 101, pp. 475–482, Jan. 2013.
%
\bibitem{Dehghan16}
S. Dehghan and N. Amjady, ``Robust transmission and energy storage expansion planning in wind farm-integrated power systems considering transmission switching,'' \emph{IEEE Trans. Sustain. Energy}, vol. 7, no. 2, pp. 765-774, Apr. 2016.
%
\bibitem{Nogales16}
A. Nogales, S. Wogrin, and E. Centeno, ``Impact of technical operational details on generation expansion in oligopolistic power markets,'' \emph{IET Gen., Trans. \& Dist.}, vol. 10, no. 9, pp. 2118-2126, Sep. 2016.
%
\bibitem{Wogrin14}
S. Wogrin, P. Dueñas, A. Delgadillo, and J. Reneses, ``New approach to model load levels in electric power systems with high renewable penetration,'' \emph{IEEE Trans. Power Syst.}, vol. 29, no. 5, pp. 2210-2218, Sep. 2014.
%
\bibitem{Merrick16}
J. H. Merrick, ``On representation of temporal variability in electricity capacity planning models,'' \emph{Energy Economics}, vol. 59, pp. 261-274, Sep. 2016.
%
\bibitem{Floudas95}
C. A. Floudas, ``Nonlinear and Mixed-Integer Optimization: Fundamentals and Applications,'' New York, NY, USA: Oxford University Press,
1995.
%
\bibitem{Cplex2016}
The ILOG CPLEX, 2016. [Online]. Available: http://www.ilog.com/ products/cplex/.
%
\bibitem{GAMS}
R.~E. Rosenthal, \emph{GAMS, A user's guide}. Washington, DC: GAMS Development Corporation, 2012.
%
\bibitem{Energinet}
Energinet - Energy Data Service, 2019. [Online]. Available: http://www.energidataservice.dk/.
%
\bibitem{RTS}
C. Grigg, P. Wong, P. Albrecht, R. Allan, M. Bhavaraju, R. Billinton, Q. Chen, C. Fong, S. Haddad, S. Kuruganty, W. Li, R. Mukerji, D. Patton, N. Rau, D. Reppen, A. Schneider, M. Shahidehpour, and C. Singh, ``The IEEE Reliability Test System-1996. A report prepared by the Reliability Test System Task Force of the Application of Probability Methods Subcommittee,'' \emph{IEEE Transactions on Power Systems}, 14(3):1010-1020, 1999.
%
\bibitem{FernandezBlanco2017}
R. Fern\'{a}ndez-Blanco, Y. Dvorkin, B. Xu, Y. Wang, and D. S. Kirschen, ``Optimal energy storage siting ans sizing: A WECC case study,'' \emph{IEEE Transactions on Sustainable Energy}, 8(2):733-743, 2017.
%

\end{thebibliography}
\end{document}